\documentclass[preprint,12pt]{elsarticle}
\usepackage[T1]{fontenc}
\usepackage[latin9]{inputenc}
\usepackage{geometry}
\usepackage{amsmath}
\usepackage{amssymb}
\usepackage{graphicx}
\usepackage{esint,color}
\usepackage{subfigure,booktabs}

\biboptions{sort&compress}

\usepackage{float}

\newtheorem{lemma}{Lemma}

\usepackage{amsfonts} 
\usepackage[english]{babel}

\begin{document}

\begin{frontmatter}

  \title{Multiscale stabilization for convection-dominated diffusion
    in heterogeneous media}

  \author{Victor M. Calo$^{1,2}$}

  \author{Eric T. Chung$^{3}$}

  \author{Yalchin Efendiev$^{1,4,5}$}

  \author{Wing Tat Leung$^{4}$}

  \address{$^{1}$ Center for Numerical Porous Media (NumPor) \\
    King Abdullah University of Science and Technology (KAUST), Thuwal
    23955-6900, Saudi Arabia.}

  \address{$^{2}$ Applied Mathematics \& Computational Science
    and Earth Science \& Engineering\\
    King Abdullah University of Science and Technology (KAUST), Thuwal
    23955-6900, Saudi Arabia.}

  \address{$^{3}$Department of Mathematics\\ The Chinese University of
    Hong Kong, Hong Kong SAR.}

  \address{$^{4}$Department of Mathematics\\
    Texas A\&M University, College Station, Texas, USA}

  \address{$^{5}$Institute for Scientific Computation (ISC) \\
    Texas A\&M University, College Station, Texas, USA}

\begin{abstract}
  We develop a Petrov-Galerkin stabilization method for multiscale
  convection-diffusion transport systems.  Existing stabilization
  techniques add a limited number of degrees of freedom in the form of
  bubble functions or a modified diffusion, which may not sufficient
  to stabilize multiscale systems. We seek a local reduced-order model
  for this kind of multiscale transport problems and thus, develop a
  systematic approach for finding reduced-order approximations of the
  solution. We start from a Petrov-Galerkin framework using optimal
  weighting functions.  We introduce an auxiliary variable to a mixed
  formulation of the problem. The auxiliary variable stands for the
  optimal weighting function. The problem reduces to finding a test
  space (a reduced dimensional space for this auxiliary variable),
  which guarantees that the error in the primal variable (representing
  the solution) is close to the projection error of the full solution
  on the reduced dimensional space that approximates the solution.  To
  find the test space, we reformulate some recent mixed Generalized
  Multiscale Finite Element Methods. We introduce snapshots and local
  spectral problems that appropriately define local weight and trial
  spaces.  In particular, we use energy minimizing snapshots and local
  spectral decompositions in the natural norm associated with the
  auxiliary variable. The resulting spectral decomposition adaptively
  identifies and builds the optimal multiscale space to stabilize the
  system.  We discuss the stability and its relation to the
  approximation property of the test space.  We design online basis
  functions, which accelerate convergence in the test space, and
  consequently, improve stability.  We present several numerical
  examples and show that one needs a few test functions to achieve an
  error similar to the projection error in the primal variable
  irrespective of the Peclet number.
\end{abstract}


\end{frontmatter}

\section{Introduction}

Existing techniques for solving multiscale problems usually seek a
reduced dimensional approximation for the solution space. Many of
these multiscale problems with high contrast require stabilization due
to the large variations in the medium properties. For example, in a
multiscale convection-dominated diffusion with a high Peclet number,
besides finding a reduced order model, one needs to stabilize the
system to avoid large errors~\cite{park2004multiscale}.  Stabilization
of multiscale methods for convection-diffusion cannot simply use a
modified diffusion and requires more sophisticated techniques.  In
this paper, we discuss a general framework for stabilization, which
combines recent developments in Generalized Multiscale Finite Element
Method (e.g.,~\cite{egh12}) and Discontinuous Petrov-Galerkin method
(e.g.,~\cite{demkowicz2014overview, niemi2011discontinuous,
  niemi2013automatically}).

We consider a convection-diffusion equation in the form
\begin{equation}
\label{eq:pde}
-\nabla\cdot(\kappa\nabla u)+b\cdot\nabla u  =f
\end{equation}
with a high Peclet number, where $\kappa$ is a diffusion tensor and
$b$ is the velocity vector~\cite{park2004multiscale,
  fannjiang1994convection}.  Both fields are characterized by
multiscale spatial features.  Many solution techniques for multiscale
problems require a construction of special basis functions on a coarse
grid~\cite{dur91, weh02, Arbogast_two_scale_04, egw10, egh12, eh09,
  g1, g2, Review, Efendiev_GKiL_12, ehg04, Chu_Hou_MathComp_10,ee03,
  calo2011note, calo2014asymptotic, GhommemJCP2013, eglmsMSDG}.  These
approaches include the Multiscale Finite Element Methods
(MsFEM)~\cite{egw10, egh12, eh09, ehg04, Ensemble, alotaibi2015global}
and Variational Multiscale Methods~\cite{hughes98,
  hughes1995multiscale, hsu2010improving, bazilevs2010isogeometric,
  masud2004multiscale, buffa2006analysis, hughes2005variational,
  codina1998comparison, bazilevs2007variational, akkerman2008role,
  hughes2007variational} among others.  In MsFEM, local multiscale
basis functions are constructed for each coarse region.  Recently, a
general framework, the Generalized Multiscale Finite Element Method
(GMsFEM), for finding a reduced approximation was
proposed~\cite{egh12, galvis2015generalized,Ensemble, eglmsMSDG,
  eglp13oversampling, calo2014multiscale, chung2014adaptive, chung2014adaptiveDG,
  randomized2014, chung2015generalizedperforated, chung2015residual,
  chung2015online}.  GMsFEM generates a reduced dimensional space on a
coarse grid that approximates the solution space by introducing local
snapshot spaces and appropriate local spectral decompositions.
However, a direct application of these approaches for
singularly-perturbed problems, such as convection-dominated diffusion,
faces difficulties due to the poor stability of these schemes.
Simplified stabilization techniques on a coarse grid are not
efficient.  Indeed, the modification of the diffusion coefficient and
similar approaches assumes the use of a few degrees of freedom
locally to stabilize the problem. These approaches do not suffice for
complex problems and one needs a systematic method to generate the
necessary test spaces.

We use the discontinuous Petrov-Galerkin (DPG) techniques
following~\cite{demkowicz2013primal, chan2014robust,
  demkowicz2013robust, demkowicz2014overview, demkowicz2012class,
  zitelli2011class} to stabilize the system.  We start with a stable
fine-scale finite element discretization that fully resolves all
scales of the underlying equation
\begin{equation}
  \label{eq:discrete}
  Au = f.
\end{equation}
The system is written in a mixed framework using an auxiliary variable
as follows
\begin{align}
  \label{eq:mixed11}
  R w + A u &= f,\\
  A^T w\ \ \qquad&=0.   \label{eq:mixed12}
\end{align}
The variable $w$ plays the role of a test function and the matrix $R$
is related to the norm in which we seek to achieve stability.  We
assume that the fine-scale system gives $w=0$, that is, it is
discretely stable.  In multiscale methods, one approximates the
solution using a reduced dimensional subspace for $u$.  More
precisely, \[u\approx\sum_{i} z^u_i \phi_i,\qquad\text{ or } \qquad
u\approx\Phi z^u.\] The resulting system also needs a reduced
dimensional test space, \[w\approx\sum_i z^w_i \psi_i\qquad\text{ or }
\qquad w\approx\Psi z^w.\] The stabilization of~\eqref{eq:discrete}
requires appropriate $\Phi$ and $\Psi$.  We discuss the design of
these spaces in the following.

Within the DPG framework, one can achieve stability by choosing
test functions $w$ with global support~\cite{barrett1984approximate,
  demkowicz1986adaptive}.  However, our goal is to design procedures
for constructing localized test spaces.  In this paper, we design a
novel test space which guarantees stability for singularly perturbed
problems such as convection-dominated diffusion in a multiscale media
with a high Peclet number.  To generate a multiscale space for $w$, we
use the recently developed theory for GMsFEM for mixed problems
\cite{chung2015mixed}.  We start by constructing a local snapshot
space which approximates the global test functions.  These snapshot
vectors are supported in coarse regions and are constructed solving
local adjoint problems in neighboring coarse elements. The snapshot
spaces are augmented with local bubble functions. The dimension of the
snapshot space is proportional to the number of fine-grid edges
(i.e., proportional to the Peclet number). To reduce the dimension of
this space, making the construction independent of the Peclet number,
we propose a set of local spectral problems.  In these local spectral
problems, we use minimum energy snapshot vectors~\cite{
  chan2015adaptive} and perform a local spectral decomposition with
respect to the $AA^T$ norm.  Our objective is to find a reduced
dimensional approximation, $w_N$, of $w$ such that $\|w-w_N\|$ is
small. We can show that the approximation property of the test space
is important to achieve the stability (cf.~\cite{dahmen2014double}).
We note that the least squares approach~\cite{bochev1998finite,
  hughes1989new, bochev2009least, fuchen16} can also be used to
achieve the stability in the natural norm. Contrary to the traditional
least squares approach, the proposed approach minimizes the residual
with some special weights related to the test functions.

We discuss how to construct online basis~\cite{ chan2015adaptive},
which uses residual information.  Online basis functions speed-up
convergence at a cost proportional to additional multiscale test
functions, which are computed by solving local problems.
In~\cite{chan2015adaptive}, we developed online basis functions for
flow equations.  One can show that by adding online basis functions,
the error reduces by a factor of $1-\Lambda_{\min}$, where $\Lambda_{\min}$ is the
smallest eigenvalue for which the corresponding eigenvector is not
included in the coarse test space. That is,
\[\|w-w_N^{\text{online}}\| \leq C
(1-\Lambda_{\min})\|w-w_N^{\text{offline}}\|,\]
where $C$ is independent of the mesh size, physical scales, and
material properties' contrast.  Thus, if we use all eigenvectors that
correspond to asymptotically small eigenvalues in the coarse test
space, it guarantees that with a few online iterations, we achieve
stability. We observe this behavior in our numerical simulations. Our
construction differs from~\cite{chan2015adaptive}. In this paper, we
design different coarse spaces for trial and test. Additionally, the
mixed formulation we present in~\eqref{eq:mixed11}-\eqref{eq:mixed12}
involves higher-order partial-derivative operators than standard mixed
forms.

Then, we present several relevant numerical examples of multiscale
transport problems.  In particular, we consider heterogeneous velocity
fields and a constant diffusion such that the resulting Peclet number
is high. We consider several types of the velocity fields.  The first
class of velocity fields we use are motivated by~\cite{
  fannjiang1994convection} and contain eddies and channels.  The
second class of velocity fields, which are motivated by porous media
applications, consist of heterogeneous channels (layers). In all
examples, we consider how the appropriate error (which is based on our
stabilization) behaves as we increase the number of test functions.
We observe that one needs several test functions per coarse degree of
freedom to achieve an error close to the projection error of the
solution of the span of the coarse degrees of freedom. Moreover, the
number of test functions does not change as we increase the Peclet
number.  By using a few test functions, we can reduce the error
achieved by standard GMsFEM by several orders of magnitude.

The paper is organized as follows. In Section 2, we present
preliminary results and notations, which include the problem setup as
well as the coarse and fine mesh descriptions. In Section 3, we
describe our proposed procedure.  Section 4 contains numerical results.
Section 5 summarizes our findings and draws conclusions.

\section{Preliminaries}

We consider the following problem
\begin{align*}
  -\nabla\cdot(\kappa\nabla u)+b\cdot\nabla u & =f,\;\text{ in \ensuremath{\Omega}}\\
  u & =0, \;\text{ on \ensuremath{\partial\Omega}}
\end{align*}
where $\kappa$ and $b$ are highly heterogeneous multiscale spatial fields
with a large ratio $\max_\Omega(b)/\min_\Omega(\kappa)$.  The weak
formulation of this problem is to find $u\in V=H^1_0(\Omega)$ such that
\[
a(u,v)=l(v), \quad\;\forall v\in V,
\]
where
\begin{align*}
  a(u,v)  &=\int_{\Omega}\kappa\nabla u\cdot\nabla v+(b\cdot\nabla u)v, \\
  l(v) &=\int_{\Omega}fv.
\end{align*}

We start with a fine-grid (resolved) discretization of the problem and
define $u_{h}$ to be the fine-grid finite element solution in the fine-grid space $V_h$, $A_{h}$
and $f_{h}$ are the stiffness matrix and the source vector on the fine
grid, that is,
\begin{align*}
  (A_{h})_{ij} & =a(\phi_{j},\phi_{i})\text{ for }\phi_{i},\phi_{j}\in V_{h}\\
  (f_{h})_{i} & =l(\phi_{i})\text{ for }\phi_{i},\phi_{j}\in V_{h}
\end{align*}
and $u_{h}=\sum\phi_{i}(u_{h})_{i}$ with $A_{h}u_{h}=f_{h}$.

We introduce an auxiliary variable (a test variable) and re-write the
system in mixed form. In particular, we consider the following
problem. Find $(u_{h}^{PG},w_{h})\in V_{h}\times V_{h}$ such that
$u_{h}^{PG}=\sum\phi_{i}(u_{h})_{i}$ and
$w_{h}=\sum\phi_{i}(w_{h})_{i}$ solve
\[
\left(\begin{array}{cc}
        A_{h}A_{h}^{T} & A_{h}\\
        A_{h}^{T} & 0
\end{array}\right)\left(\begin{array}{c}
                          w_{h}\\
                          u_{h}^{PG}
\end{array}\right)=\left(\begin{array}{c}
                           f_{h}\\
                           0
\end{array}\right).
\]
Since $det(A_{h})\neq0$, we have $u_{h}^{PG}=u_{h}$ and $w_{h}=0$.
Therefore, these two problems have the same solution. Our objective is
to find a reduced dimensional coarse approximation for $w_{h}$, which
can guarantee that the corresponding $u_{h}^{PG}$ is a good
approximation to $u_{h}$.

\begin{figure}[htb]
  \centering
  \includegraphics[width=0.65 \textwidth]{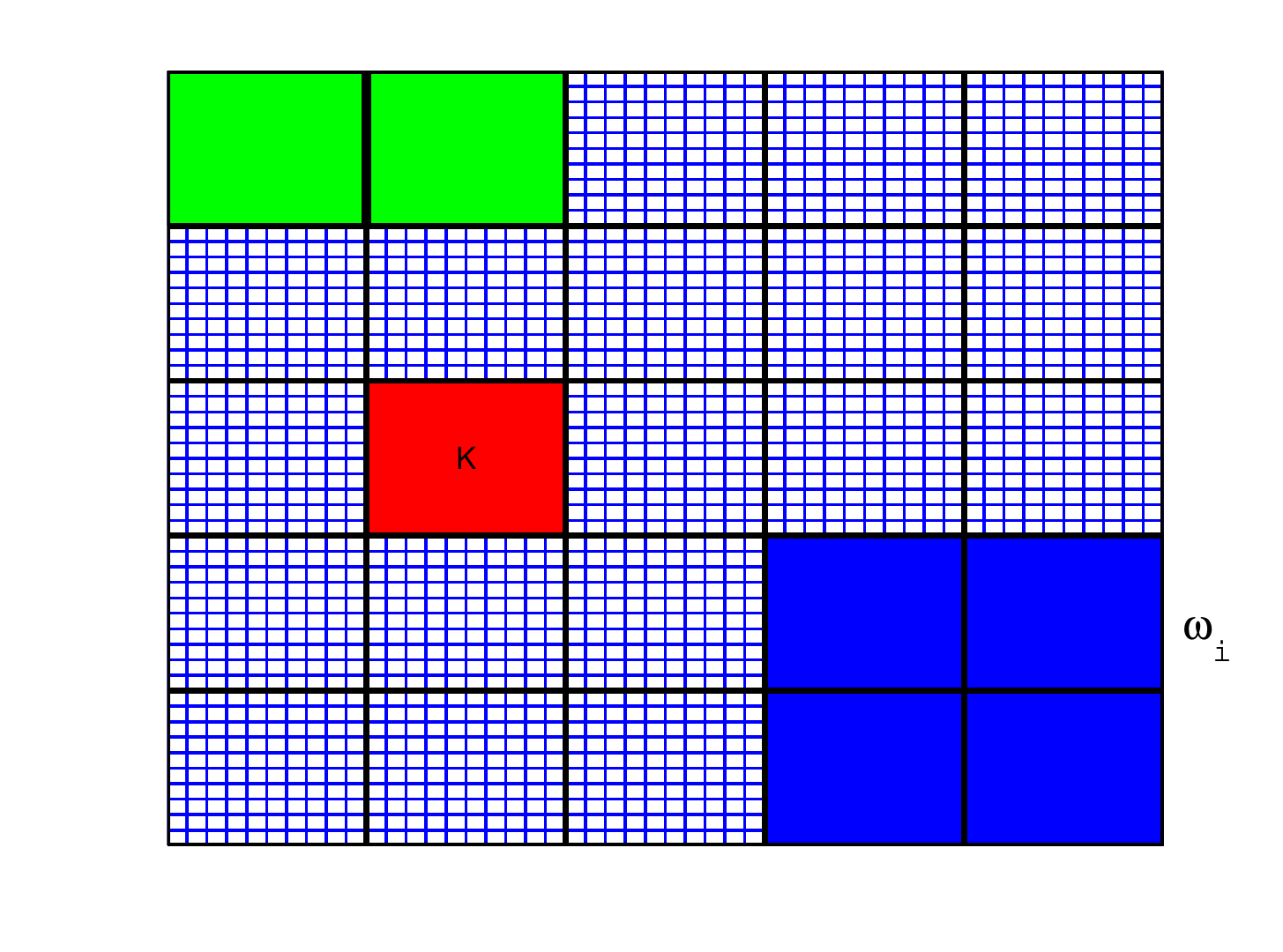}
  \caption{Illustration of coarse neighborhoods and
    elements. Red designates a coarse element. Green designates two
    neighboring elements that share a common face (used to construct
    test functions). Blue designates the union of all coarse elements
    that share a common vertex (used to construct trial
    functions).}
  \label{schematic}
\end{figure}

\subsection{Coarse-grid description}

Next, we introduce some notation.  We use $\mathcal{T}^H$ to denote a
conforming partition of the computational domain $D$.  The set
$\mathcal{T}^H$ is called the coarse grid and the elements of
$\mathcal{T}^H$ are called coarse elements.  Moreover, $H>0$ is the
coarse mesh size.  We only  consider rectangular coarse
elements to simplify the  discussion and illustrations.  The
methodology presented can be easily extended to coarse elements with
more general geometries.  Let $N$ be the number of nodes in the coarse
grid $\mathcal{T}^H$, and let $\{x_i \, | \, 1\leq i\leq N\}$ be the
set of nodes in the coarse grid (or coarse nodes for short).  For each
coarse node $x_i$, we define a coarse neighborhood $\omega_i$ by
\begin{equation} \label{neighborhood} \omega_i=\bigcup\{
  K_j\in\mathcal{T}^H; ~~~ x_i\in \overline{K}_j\}.
\end{equation}
That is,
$\omega_i$ is the union of all coarse elements $K_j\in\mathcal{T}^H$
having the coarse node $x_i$ (blue region in
Figure~\ref{schematic}). We use two neighboring elements sharing a
common face to construct the test functions. An example of this region
is depicted in green in Figure~\ref{schematic}.

\section{Generalized Multiscale Finite Element Method for
  Petrov-Galerkin Approximations}

In this section, we discuss the construction of the multiscale basis
functions for the trial space $V$ and the test space $W$.  In
particular, we show that one needs a good approximation for $w_h$ in
order to achieve discrete stability.  We start by introducing some
notation and formulating the multiscale Petrov-Galerkin framework we
solve.  We introduce the snapshot space and then the local
spectral decomposition used to construct the multiscale basis
functions.

To simplify notations, we let
\begin{equation*}
Au := -\nabla\cdot(\kappa\nabla u)+b\cdot\nabla u,
\end{equation*}
and
\begin{equation*}
A^*u := -\nabla\cdot(\kappa\nabla u)- \nabla\cdot (b \,u).
\end{equation*}

\begin{figure}[htb!]
  \centering
  \subfigure[Trial snapshot basis]{\label{fig:illtrialLeft}
    \includegraphics[width = 0.475\textwidth]{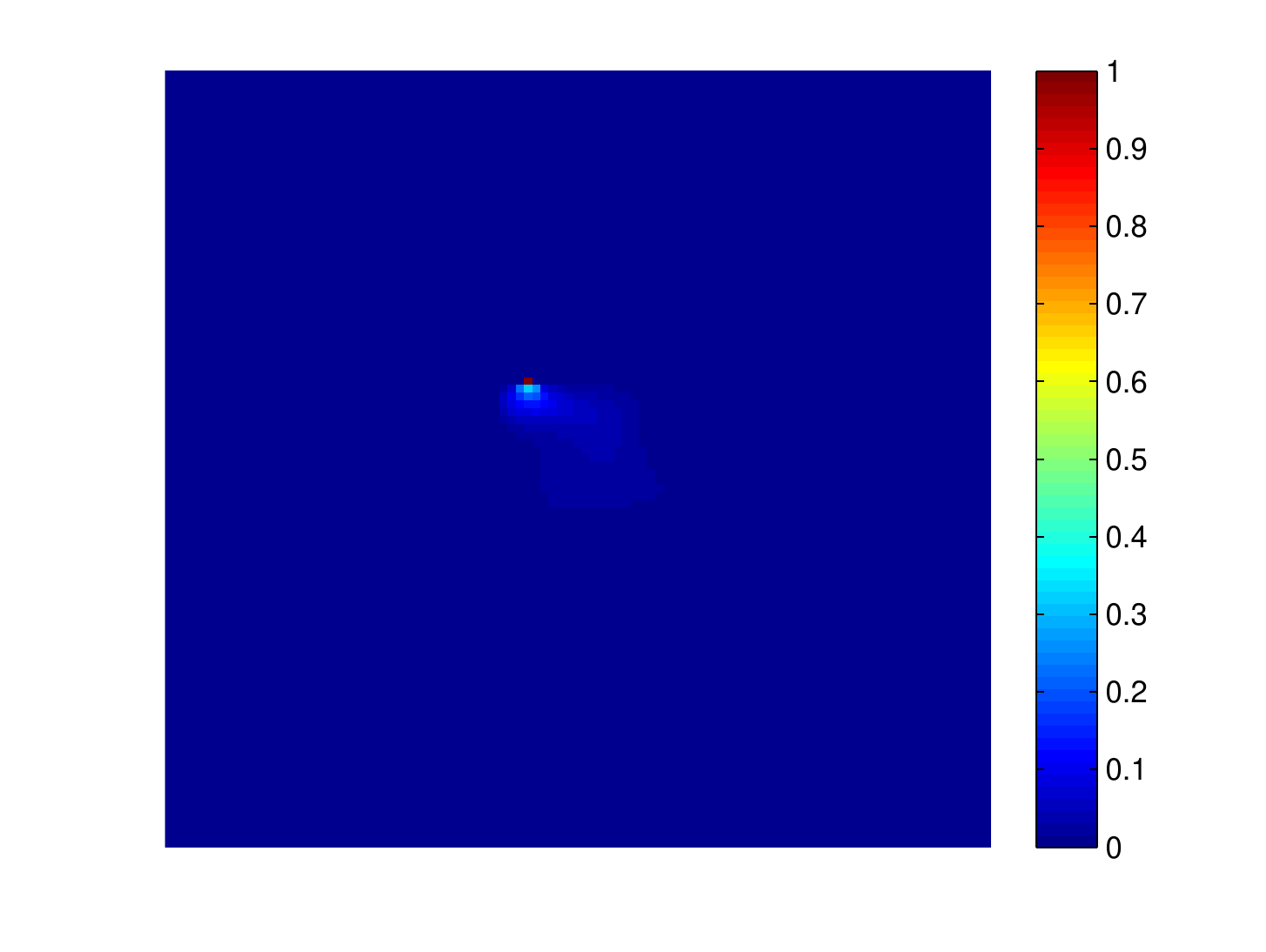}
   }
  \subfigure[ Offline trial basis]{\label{fig:illtrialRight}
     \includegraphics[width = 0.475\textwidth]{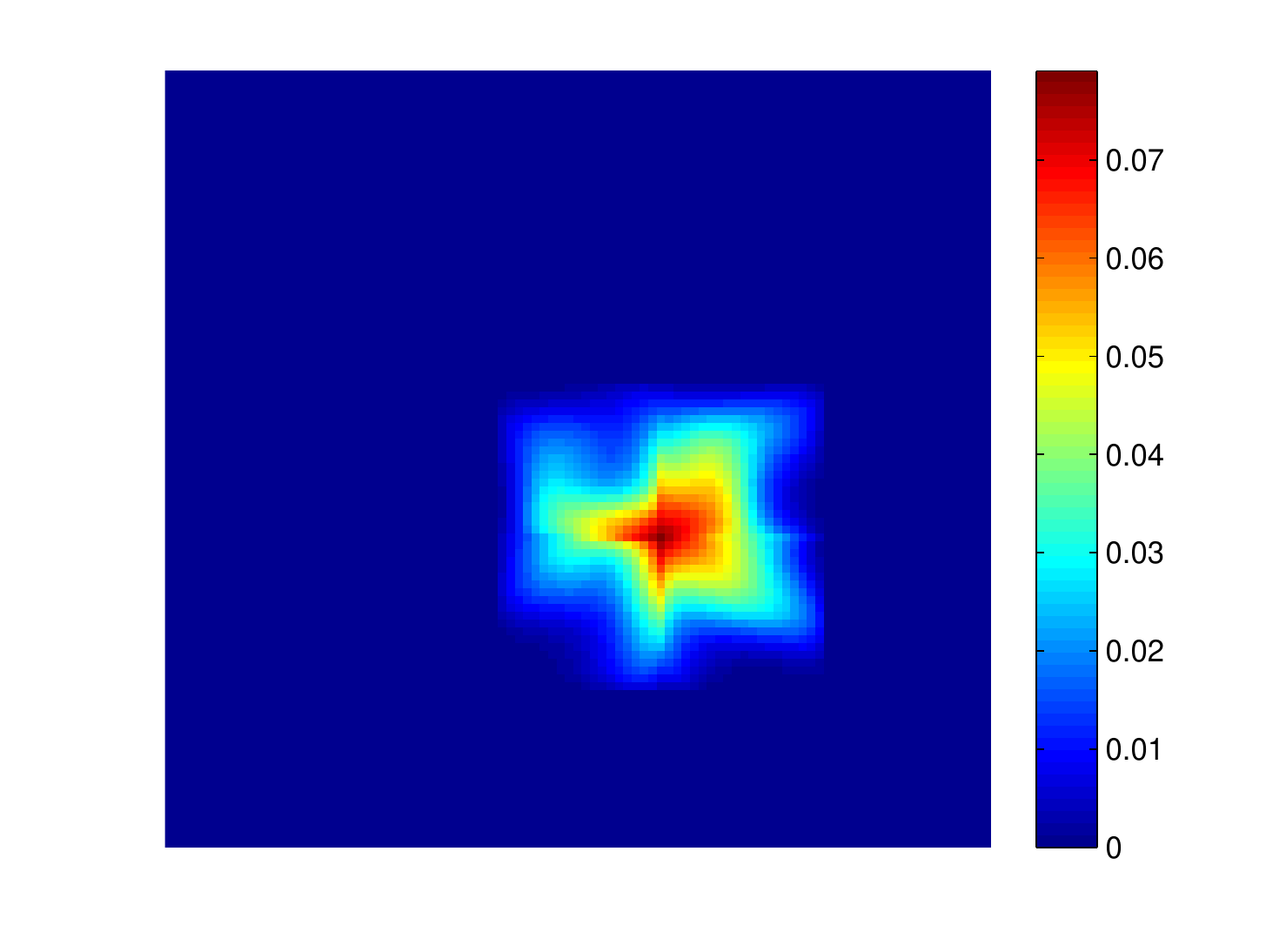}
  }
  \caption{A trial snapshot basis and the resulting offline trial
    basis for a given coarse block.}
  \label{fig:illtrial}
\end{figure}

\subsection{Construction of the multiscale trial space}

\subsubsection{Snapshot space}
We solve a local problem with specifically designed boundary
conditions to construct the snapshot basis functions.
For each coarse neighborhood $\omega_l$,
 we define a set of snapshot functions
$\phi_{i,l}^{snap}$ such that
\begin{align*}
A(\phi_{i,l}^{snap}) :=
  -\nabla\cdot(\kappa\nabla\phi_{i,l}^{snap})+b\cdot\nabla\phi_{i,l}^{snap}&=0,&&\text{in }\omega_{l},\\
  \phi_{i,l}^{snap}|_{\omega_{l}}(x_{j})&=\delta_{ij},&&\text{on
  }\partial\omega_{l},
\end{align*}
where $\delta_{ij}$ is the discrete delta function defined on $\partial\omega_l$
with respect to the fine grid.
The local snapshot space for the trial space is defined by
$V^{snap}(\omega_{l})=\text{span}\{\phi_{i,l}^{snap}\}$.  The snapshot
functions and multiscale basis functions (offline space) are defined
in the union of coarse elements that share a common vertex
(Figure~\ref{schematic} shows a schematic representation of the grid while
Figure~\ref{fig:illtrial} shows a solution snapshot and a multiscale
basis function).
We use $\Phi_l^T$ to denote the change of basis matrix
from the fine-grid space $V_h(\omega_l)$ to $V^{snap}(\omega_l)$.
Here $V_h(\omega_l)$ is the restriction of $V_h$ in $\omega_l$.

\subsubsection{Eigenproblem}

To construct the offline trial space, we solve the following
eigenproblem
\[
(A_{snap}^{\omega_{l}})^{T}A_{snap}^{\omega_{l}}v_{j}=\lambda_j
M_{snap}^{\omega_{l}}v_{j},
\]
where
\begin{align*}
  A_{snap}^{\omega_{l}} & =\Phi_{l}^T A_{h}^{\omega_{l}}\Phi_{l}\\
  M_{snap}^{\omega_{l}} & =\Phi_{l}^T M_{h}^{\omega_{l}}\Phi_{l}
\end{align*}
and $(\lambda_j,v_j)$ is the $j$-th eigen-pair.
In the above definition, $A_{h}^{\omega_{l}}$ and $M_{h}^{\omega_{l}}$
are the restrictions of the fine-scale stiffness matrix $A_h$ and the fine-scale mass matrix $M_h$
in $\omega_l$.
We order the eigenvalues in increasing order
and we use the first $m$ eigenfunctions as the offline trial basis functions.
Specifically, we define
$\xi_{l,j}=\Phi_{l}v_{j}$ and
$V_{off}=\text{span}\{\chi_{l}\xi_{l,j}|1\leq j\leq m, 1\leq l \leq N\}$, where
$\{\chi_{l}\}$ is the partition of unity.

\begin{figure}[htb!]
  \centering
  \subfigure[ Test snapshot basis for $W_{3}$]{\label{fig:illtestLeft}
    \includegraphics[width = 0.475\textwidth]{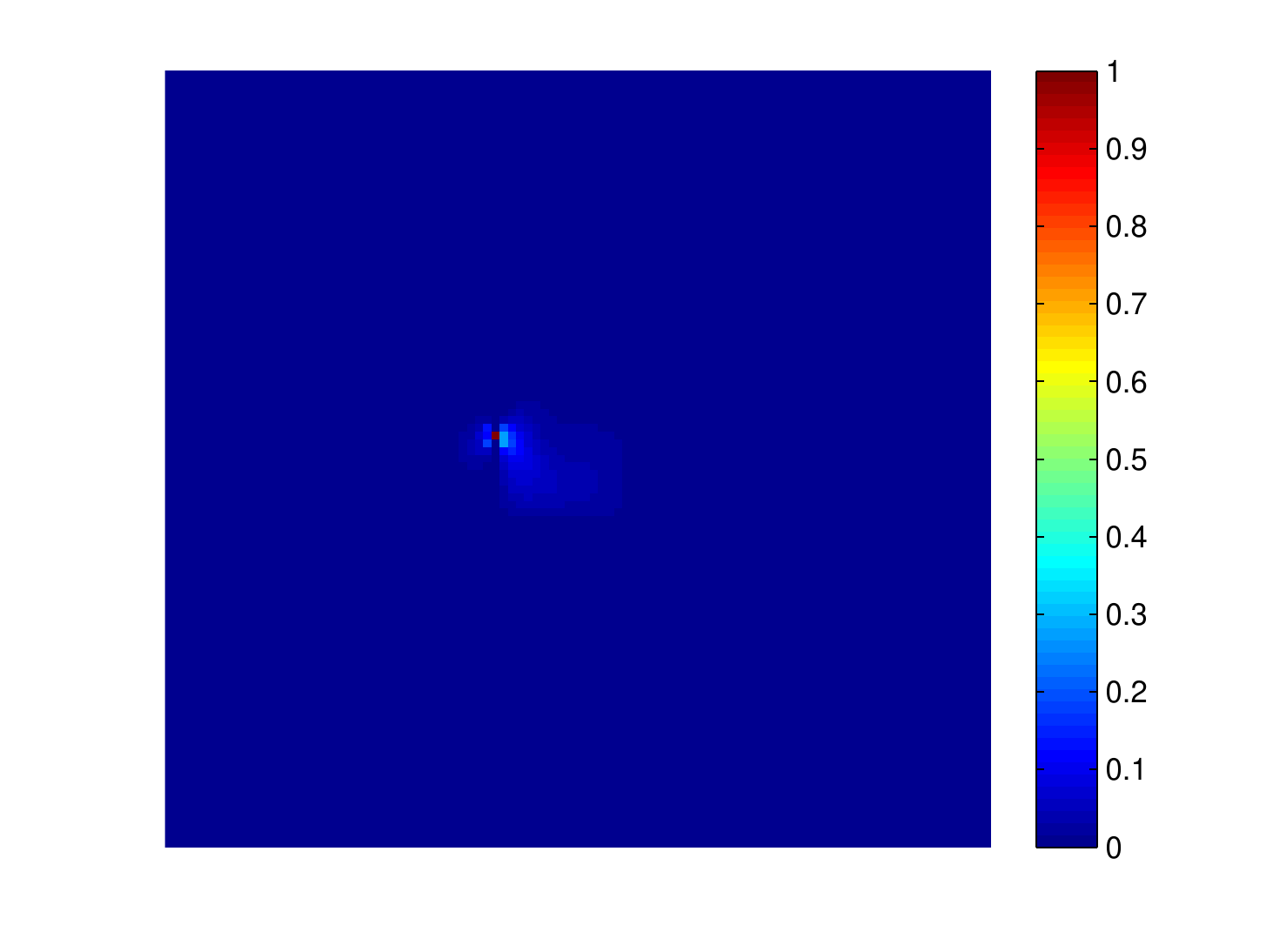}
   }
  \subfigure[ Offline test basis for $W_{3}$]{\label{fig:illtestRight}
     \includegraphics[width = 0.475\textwidth]{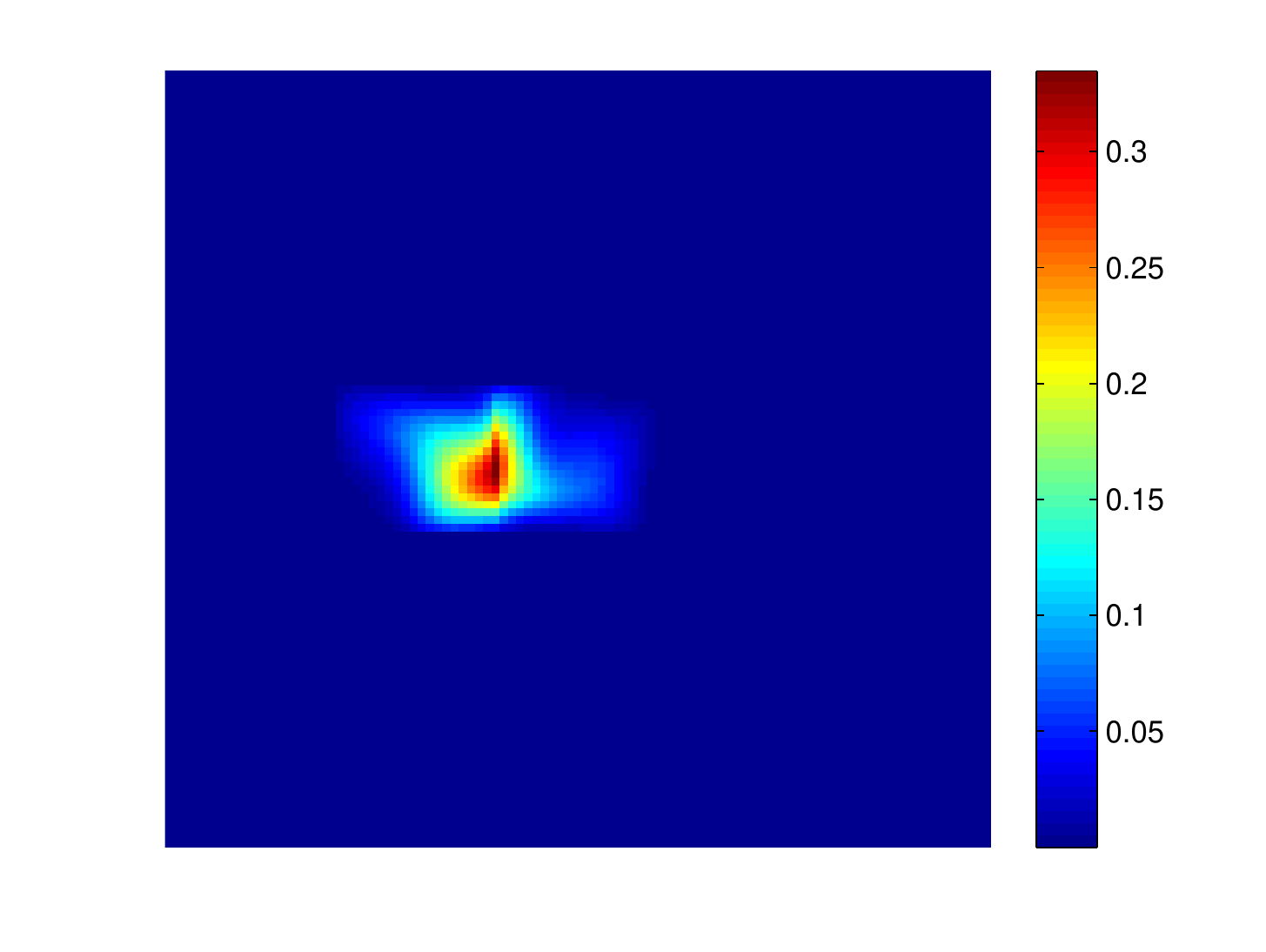}
  }
  \caption{A test snapshot basis and the resulting offline test
    basis for a given coarse block of $W_{3}$.}
  \label{fig:illtest}
\end{figure}

\subsection{Construction of the multiscale test space}

\subsubsection{Snapshot space}

The snapshot space for the test space consists of three components, and is denoted as $W_{snap}=W_{snap}^{1}+W_{snap}^{2}+W_{snap}^{3}$.
Next, we will give the constructions for $W_{snap}^1, W_{snap}^2$ and $W_{snap}^3$.
In each coarse block $K_{k}$, we define $W_{snap}^{1}(K_{k})$ as
\[
W_{snap}^{1}(K_{k}):=\{\psi^{snap}\in
V_{h,0}(K_{k}) \; | \; A^*(\psi^{snap})=\xi_{l,j}\;\text{
  in }K_{k}\text{ for some }\xi_{l,j}\in V_{off}\}
\]
where $V_{h}(K_k)$ is the restriction of $V_h$ in $K_k$
and $V_{h,0}(K_k)$ is the subspace of $V_h(K_k)$ containing functions that vanish on $\partial K_k$.
The space $W_{snap}^{1}(K_{k})$ contains functions that are solution
of the adjoint problem on $K_k$ with a source term $\xi_{l,j}$
and zero Dirichlet boundary condition.
The space $W_{snap}^{1}$, is defined as
$W_{snap}^{1}=\oplus_{k}W_{snap}^{1}(K_{k})$.
The space $W_{snap}^{1}$ is considered as the space of multiscale bubble functions.
We remark that we obtain perfect test functions (with perfect stabilization) if the above local problems
are solved on the whole domain.

The second space $W_{snap}^{2}$ is defined as follows.
For each coarse block $K_k$, we define
\begin{equation*}
W_{snap}^{2}(K_k):=\{\psi^{snap} \in V_h(K_k) \;|\; A^*(\psi^{snap})=0\;\text{
  in }K_{k}\text{ and }\psi^{snap}\text{ is linear on } E \in \partial
K_{k}\}.
\end{equation*}
The space $W_{snap}^{2}$ is defined by $W_{snap}^{2}=\oplus_{k}W_{snap}^{2}(K_{k}) \cap C^0(\Omega)$.
Note that this space is similar to the classical multiscale finite element space.
%\[
%W_{snap}^{2}:=\{\psi^{snap}|-\nabla\cdot(\kappa\nabla\psi^{snap})-\nabla\cdot(b\psi^{snap})=0\;\text{
  %in }K_{k}\text{ and }\psi^{snap}\text{ is linear on }\partial
%K_{k}\text{ for all }K_{k}\}.
%\]

Finally, we give the definition for $W_{snap}^{3}$.
For each coarse edge $E_{k}$,
we define $K(E_k)$ as the set of all coarse blocks having the edge $E_k$.
Then, we find $\psi_{i,k}^{snap} \in V_h(K(E_k))$ such that
\begin{align*}
  A^*\psi_{i,k}^{snap} := -\nabla\cdot(\kappa\nabla\psi_{i,k}^{snap})-\nabla\cdot(b\psi_{i,k}^{snap}) & =0\;\text{ in each } K \in K(E_{k}), \\
  \psi_{i,k}^{snap}|_{E_{k}}(x_{j}) & =\delta^0_{ij}\text{ for all }x_{j}\in E_{k},\\
  \psi_{i,k}^{snap}|_{\partial K(E_{k})\backslash E_{k}} & =0.
\end{align*}
In the above system, $\delta^0_{ij}$ is the discrete delta function defined on $E_k$
with respect to the fine mesh and is zero on the boundary of $E_k$.
Then we define $W_{snap}^{3}(E_{k})=\text{span}\{\psi_{i,k}^{snap}\}$, and
$W_{snap}^{3}=\oplus_{k}W_{snap}^{3}(E_{k})$.  (Figure~\ref{schematic}
illustrates a grid and Figure~\ref{fig:illtest} illustrates a snapshot
solution and a multiscale basis function).

We remark that $W_{snap} = V_h$.

\begin{lemma}
  For each $u\in V_\text{{off}}$, there exists a test function
  $\phi\in W_\text{{snap}}$ such that
\begin{equation}
  a(v,\phi) = (u,v)_{l^2}, \; \forall v\in V_h.
\label{eq:lemma1}
\end{equation}
\end{lemma}
{\bf Proof:}
Given $u\in V_\text{{off}}$, we assume that $\phi \in V_h $ satisfies
\eqref{eq:lemma1}.  For each $K\in\mathcal{T}^{H}$, we define
$\phi_{\text{snap}}^{(1)}\in W_{\text{snap}}^{1}$ satisfying
\[
a(v,\phi_{\text{snap}}^{(1)})=(u,v)_{l^2},\quad \forall v\in V_{h,0}(K),
\; \forall K\in\mathcal{T}^{H}.
\]
Next, we define
$\phi_{\text{snap}}^{(2)}\in W_{\text{snap}}^{2}+W_{\text{snap}}^{3}$
such that
\[
\phi_{\text{snap}}^{(2)}=\phi\text{ on \ensuremath{\partial}K}, \;
\forall K\in\mathcal{T}^{H}.
\]

Then, we have
\[
a(v,\phi-\phi_{\text{snap}}^{(1)}-\phi_{\text{snap}}^{(2)})=0, \;
\forall v\in V_{h,0}(K), \; \forall K\in\mathcal{T}^{H}
\]
with
$\phi-\phi_{\text{snap}}^{(1)}-\phi_{\text{snap}}^{(2)}=0\text{ on
  \ensuremath{\partial}K}, \; \forall K\in\mathcal{T}^{H}.$
Therefore, $\phi-\phi_{\text{snap}}^{(1)}-\phi_{\text{snap}}^{(2)}$ is
a solution of the adjoint problem with zero Dirichlet boundary
condition and zero source term in all $K\in\mathcal{T}^{H}$. Thus,
$\phi=\phi_{\text{snap}}^{(1)}+\phi_{\text{snap}}^{(2)}\in
W_{\text{snap}}.$

\subsubsection{Eigenproblem}
\label{sec:TestEigenproblem}
Among the three parts of the test space, the dimension of
$W_{snap}^{3}$ is proportional to the number of fine-grid blocks and
thus proportional to the Peclet number of the problem.  Consequently,
our objective is to reduce the degrees of freedom associated with
$W_{snap}^{3}$.  Both the dimensions of $W_{snap}^{1}$ and
$W_{snap}^{2}$ are proportional to the number of coarse grid degrees
of freedom.  We consider two different eigenvalue problems to
construct the offline test space for $W_{snap}^{3}$.

{\bf The first eigenvalue problem for $W_{snap}^{3}(E_{k})$:} In this
eigenvalue problem, we will use the edge values of the snapshot
solutions.
\[
\int_{K(E_{k})}(A^{T}v)(A^{T}\psi_{j})=\lambda\int_{E_{k}}v\psi
\]
The eigenvalues go to $\infty$ as we refine the fine mesh.

{\bf The second eigenvalue problem for $W_{snap}^{3}(E_{k})$:} This
eigenvalue problem is motivated by~\cite{chan2015adaptive}, where we
construct minimum energy snapshot solutions and perform a local
spectral decomposition using the same norms. More precisely,
\[
\int_{K(E_{k})}(A^{T}\tilde{v})(A^{T}\tilde{\psi_{j}})=\lambda\int_{K(E_{k})}(A^{T}v)(A^{T}\psi_{j})
\]
where $\tilde{\psi}=\text{argmin}_{\tilde{\psi}\in\{v\in
  W_{snap}^{3}|v|_{E_{k}}=\psi|_{E_{k}}\}}\{\int_{K(E_{k})}(A^{T}\tilde{\psi})(A^{T}\tilde{\psi})\}$.
In this case, the eigenvalues are always smaller than $1$.

We will arrange the eigenvalues of the above spectral problems in increasing order,
and choose the first $L_k$ eigenfunctions as the offline test basis functions. The
span of these basis functions is denoted as $W_{off}^3$.
The final test space $W_{off}$ is defined by $W_{snap}^1 \oplus W_{snap}^2 \oplus W_{off}^3$.

\subsection{Global coupling}

We can use the above trial and test spaces to get a reduced system for the multiscale solution.
In particular,
the multiscale solution is computed by solving
\begin{equation}
\label{eq:global}
\left(\begin{array}{cc}
        \Theta^{T}A_{h}A_{h}^{T}\Theta & \Theta^{T}A_{h}\Xi\\
        \Xi^{T}A_{h}^T\Theta & 0
      \end{array}\right)\left(\begin{array}{c}
                                w_{ms}\\
                                u_{ms}^{PG}
                              \end{array}\right)=\left(\begin{array}{c}
                                                         \Theta^T f_{h}\\
                                                         0
                                                       \end{array}\right).
\end{equation}
The columns of $\Theta$ consist of the computed multiscale test
functions while the columns of $\Xi$ consist of the computed
multiscale trial functions.

\subsection{Discussion}

Next, we discuss the approximation properties for the test space and
how the definition of this space affects the discrete stability of the
resulting method.  To simplify the discussion, we introduce some
notation. Let $N$ and $M$ be the dimensions for the test and trial
spaces, respectively. Thus, we can write
\begin{equation}
\label{eq:not1}
\left(\begin{array}{cc}
        \Theta_N^{T}A_{h}A_{h}^{T}\Theta_N & \Theta_N^{T}A_{h}\Xi_M\\
        \Xi_M^{T}A^T_{h}\Theta_N & 0
\end{array}\right)\left(\begin{array}{c}
                          w_{N,M}\\
                          u_{N,M}^{PG}
\end{array}\right)=\left(\begin{array}{c}
                           \Theta_N^{T} f_{h}\\
                           0
\end{array}\right).
\end{equation}
For simplicity, we denote by $\Theta_\infty$ the snapshot matrix that
contains all snapshot vectors in the test space and similarly for the
trial space $\Xi_\infty$. Therefore, the following statements are
true:

\begin{itemize}
\item
$w_{\infty,\infty}=0$.

\item $u_{\infty,M}$ is a projection of $u_{\infty,\infty}$ onto
  $\Xi_M$
\[
u_{\infty,M}=\Pi_{\Xi_M} u_{\infty,\infty}.
\]

\item Our objective is to find the smallest possible $N$ and $M_0$,
  such that
  $\|u^{PG}_{N,M}-u^{PG}_{\infty,\infty}\|\preceq
  \|u^{PG}_{\infty,\infty}-u^{PG}_{\infty,M}\|$
  for any $M$, when $M>M_0$.

\item The inf-sup condition for our discrete saddle-point problem can
  be written as
\begin{equation}
  \sup_\Theta \frac{\Theta_N^{T}A_{h}\Xi_M}{(\Theta_N^{T}A_{h}A_{h}^{T}\Theta_N)^{1/2}}\geq C_{infsup} (\Xi_M^T M \Xi_M)^{1/2}.
\end{equation}
The inf-sup condition implies that
\[
C_{infsup}=\inf_{u=\Xi_M q} \frac{\|\Pi_{A^T \Theta_N}
  (u)\|_{l^2}}{\|u\|_{l^2}}.
\]
Next, we take $u=A^T z$. Then, the projection of $u$ onto
$A^T \Theta_N$ is
\begin{align}
  \Pi_{A^T \Theta_N} (u)&=A^T \Theta_N ((A^T \Theta_N)^T A^T \Theta_N)^{-1} (A^T \Theta_N)^T A^T z\nonumber\\
  &=A^T \Theta_N (\Theta_N^T A A^T \Theta_N)^{-1} \Theta_N^T A A^T z.
\label{eq:choicez}
\end{align}
We define the $\Theta$ projection in the $B$ norm to be
\[
\Pi_{\Theta, B}(z)=\Theta (\Theta^T B \Theta)^{-1} \Theta^T B z.
\]
Thus,
\begin{align}
  \|\Pi_{A^T \Theta_N} (u)\|^2&= z^T A A^T \Theta_N (\Theta_N^T A A^T
  \Theta_N)^{-1} \Theta_N^T A
  A^T \Theta_N (\Theta_N^T A A^T \Theta_N)^{-1} \Theta_N^T A A^T z\nonumber\\
  &=z^T A A^T \Theta_N (\Theta_N^T A A^T \Theta_N)^{-1}  \Theta_N^T A A^T z\nonumber\\
  &=\|\Pi_{\Theta_N,AA^T}(z)\|_{AA^T}^2.
\end{align}
Also,
\[
\|u\|_{l^2}=\|z\|_{AA^T}.
\]
Thus,
\begin{equation}
\label{eq:form_infsup}
C_{infsup}=\inf_{u=A^Tz, u = \Xi_M q} \frac{\|\Pi_{\Theta_N,AA^T}(z)\|_{AA^T}} {\|z\|_{AA^T} }.
\end{equation}
 If the inf-sup is satisfied, then we have
\begin{align}
  \|w_{N,M}-w_{\infty,\infty}\| +
  \|u^{PG}_{N,M}-u^{PG}_{\infty,\infty}\|&\preceq
  \|\widehat{w}_{N,M}-w_{\infty,\infty}\| + \|\widehat{u^{PG}}_{N,M}-u^{PG}_{\infty,\infty}\|\nonumber\\
  &=0 + \|{u^{PG}}_{\infty,M}-u^{PG}_{\infty,\infty}\|.
\end{align}
From here, we have
\begin{equation}
\begin{split}
  \|w_{\infty,M}-0\| +
  \|u^{PG}_{\infty,M}-u^{PG}_{\infty,\infty}\|\preceq
  \|\widehat{u^{PG}}_{\infty,M}-u^{PG}_{\infty,\infty}\|.
\end{split}
\end{equation}
Because $N=\infty$, ${u^{PG}}_{\infty,M}=\widehat{u^{PG}}_{\infty,M}$,
we get
\[
\|w_{\infty,M}\| \preceq
\|\widehat{u^{PG}}_{\infty,M}-u^{PG}_{\infty,\infty}\|.
\]

\item The discrete inf-sup condition can be shown if for any $z$
  (e.g., $z=A^{-T} \Xi_M q$), there exists $z_0$ in the space spanned
  by $\Theta_N$ (i.e., $z_0=\Theta_N z_r$), such that
\begin{equation}
  \|z-z_0\|_{AA^T} \leq \delta \|z\|_{AA^T},
\end{equation}
for some $\delta<1$. In multiscale methods (in particular, in our
works~\cite{chan2015adaptive, chung2015mixed}), we reduce the error in
$\|z-z_0\|_{AA^T}$ by selecting appropriate multiscale spaces (as
those used herein).  In addition, this procedure can be done
adaptively.  Thus, by selecting a sufficient number of multiscale
basis functions, we can reduce the error $\|z-z_0\|_{AA^T}$ and can
achieve the stability sought.  We do not have rigorous error
estimates, but study this problem numerically. We emphasize that we
need good approximation properties in the test space (as in~\cite{
  dahmen2014double}), which is due to the primal formulation and the
choice of $z$ in~\eqref{eq:choicez}.

\end{itemize}

\subsection{Online test basis construction (residual-driven
  correction)}

One can use residual information to construct online basis
functions. Online basis functions use global information and thus
accelerate the convergence.  In~\cite{ chan2015adaptive}, we discuss
the online basis construction for flow equations using a mixed
formulation.  We use the local residual to construct an online basis
function locally in each non-overlapping coarse grid region
$\omega_{i}.$

The offline solution in the fine-scale test space
$(w_{\infty,M},u_{\infty,M})\in V_{h}\times V_{off}$ satisfies
\begin{equation}
  \left(\begin{array}{cc}
      A_{h}A_{h}^{T} &  A_{h} \Xi_M\\
      \Xi^{T}_M A_{h}^{T} & 0
\end{array}\right)\left(\begin{array}{c}
  w_{\infty,M}\\
  u_{\infty,M}
\end{array}\right)=\left(\begin{array}{c}
   f_{h}\\
  0
\end{array}\right)
\end{equation}
and the multiscale solution $(w_{N,M},u_{N,M})\in W_{off}\times
V_{off}$ satisfies
\begin{equation}
  \left(\begin{array}{cc}
      \Theta^{T}_N A_{h}A_{h}^{T}\Theta_N & \Theta^{T}_N A_{h}\Xi_M\\
      \Xi^{T}_M A_{h}^{T}\Theta_N & 0
\end{array}\right)\left(\begin{array}{c}
  w_{N,M}\\
  u_{N,M}
\end{array}\right)=\left(\begin{array}{c}
  \Theta^{T}_N f_{h}\\
  0
\end{array}\right).
\label{eq:coarseeq}
\end{equation}

The above motivates the following
local residual operator $R_i$, which is defined as
$R_{i}:V_{h}(\omega_{i})\rightarrow\mathbb{R}$ is defined by
\[
R_{i}(v)=v^{T}\Big((A_{h}A_{h}^{T})^{(i)} \Theta_N w_{N,M}+(A_{h})^{(i)} \Xi_M u_{N,M}-f_{h}\Big)
\]
and the local residual norm, $\|R_{i}\|$ is defined by
\[
\|R_{i}\|=\sup_{v\in
  V_{h}(\omega_{i})}\cfrac{|R_{i}(v)|}{\sqrt{v^{T}(A_{h}A_{h}^{T})^{(i)}v}} \, ,
\]
where $(A_{h}A_{h}^{T})^{(i)}$ and $A_{h}^{(i)}$ are local sub-matrices
of $A_{h}A_{h}^{T}$ and $A_{h}$ which correspond to the coarse grid
subdomain $\omega_{i}$.  Next, we use the local residual to construct
the local test basis, $\phi_{on}^{(i)}\in V_{h}(\omega_{i})$ such that
\[
v^{T}(A_{h}A_{h}^{T})^{(i)}\phi_{on}^{(i)}=R_{i}(v), \quad \forall v\in
V_{h}(\omega_{i}).
\]
In~\cite{chan2015adaptive}, we show that if online basis functions are
constructed using the second eigenvalue problem, then the error will
decrease at a rate $(1-\min_E\Lambda_{\min}^E)$,
where $\Lambda_{\min}^E$ is the minimum
of the eigenvalues of the spectral problem defined on $W_{snap}^3(E)$
corresponding to eigenfunctions not chosen as basis.
Consequently, using
online basis functions, we can achieve the discrete inf-sup stability
in one iteration provided $\min_E\Lambda_{\min}^E>0$.

\subsubsection{Online test basis enrichment algorithm}
First, we choose an offline trial space, $V_{\text{off}}$ and an
initial offline test space, $W^{(1)}_{\text{off}}$, by fixing the number
of basis functions for each coarse neighborhood.  Next, we construct a
sequence of online test spaces $W^{(m)}_{\text{off}}$ and compute the multiscale
solution $(w^{(m)}_{\text{ms}},u^{(m)}_{\text{ms}})$ by solving equation
\eqref{eq:coarseeq}.  The test space is constructed iteratively for $
m = 1,2,3,\dots,$ by the following algorithm:
\begin{enumerate}
\item[Step 1:] Find the multiscale solution in the current
  space. Solve for $(w^{(m)}_{\text{ms}},u^{(m)}_{\text{ms}})\in
  W^{(m)}_{\text{off}}\times V_{\text{off}}$ such that
\[
\left(\begin{array}{cc}
    (\Theta^{(m)}_{\text{off}})^{T}A_{h}A_{h}^{T} \Theta^{(m)}_{\text{off}} & (\Theta^{(m)}_{\text{off}})^{T}A_{h} \Pi^{(m)}_{\text{off}} \\
   ( \Xi^{(m)}_{\text{off}})^{T}A_{h}^{T} \Theta^{(m)}_{\text{off}} & 0
\end{array}\right)\left(\begin{array}{c}
  w^{(m)}_{\text{ms}}\\
  u^{(m)}_{\text{ms}}
\end{array}\right)=\left(\begin{array}{c}
  (\Theta^{(m)}_{\text{off}})^{T}f_{h}\\
  0
\end{array}\right).
\]

\item[Step 2:] For each coarse region $\omega_i$, compute the online
  basis, $\phi^{(i)}_{\text{on}}\in V_h(\omega_i)$ such that
\[
v^{T}(A_{h}A_{h}^{T})^{(i)}\phi_{on}^{(i)}=R_{i}(v), \quad \forall v\in
V_{h}(\omega_{i}).
\]

\item[Step 3:] Enrich the test space by setting
\[
W^{(m+1)}_{\text{off}} = W^{(m)}_{\text{off}} + \text{span}\{
\phi^{(i)}_{\text{on}}\}.
\]

\end{enumerate}

We remark that in each iteration, we perform the above procedure on non-overlapping
coarse neighborhoods, see \cite{chan2015adaptive}.

\section{Numerical Results}

\begin{figure}[!ht]
  \centering
  \subfigure[ $\alpha=2$]{\label{fig:solns_ex1L}
    \includegraphics[width = 0.475\textwidth]{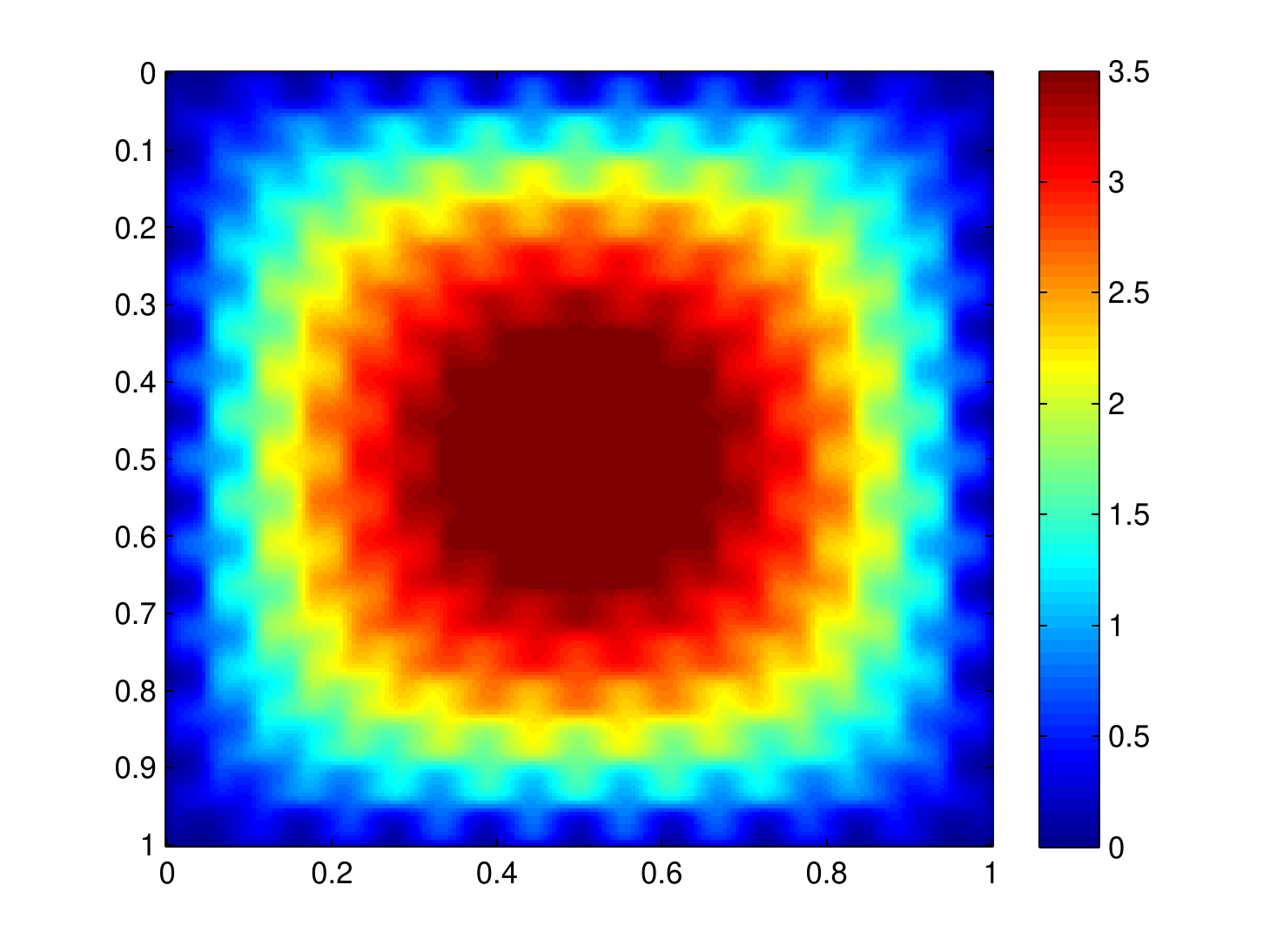}
   }
  \subfigure[ $\alpha=4$]{\label{fig:solns_ex1R}
     \includegraphics[width = 0.475\textwidth]{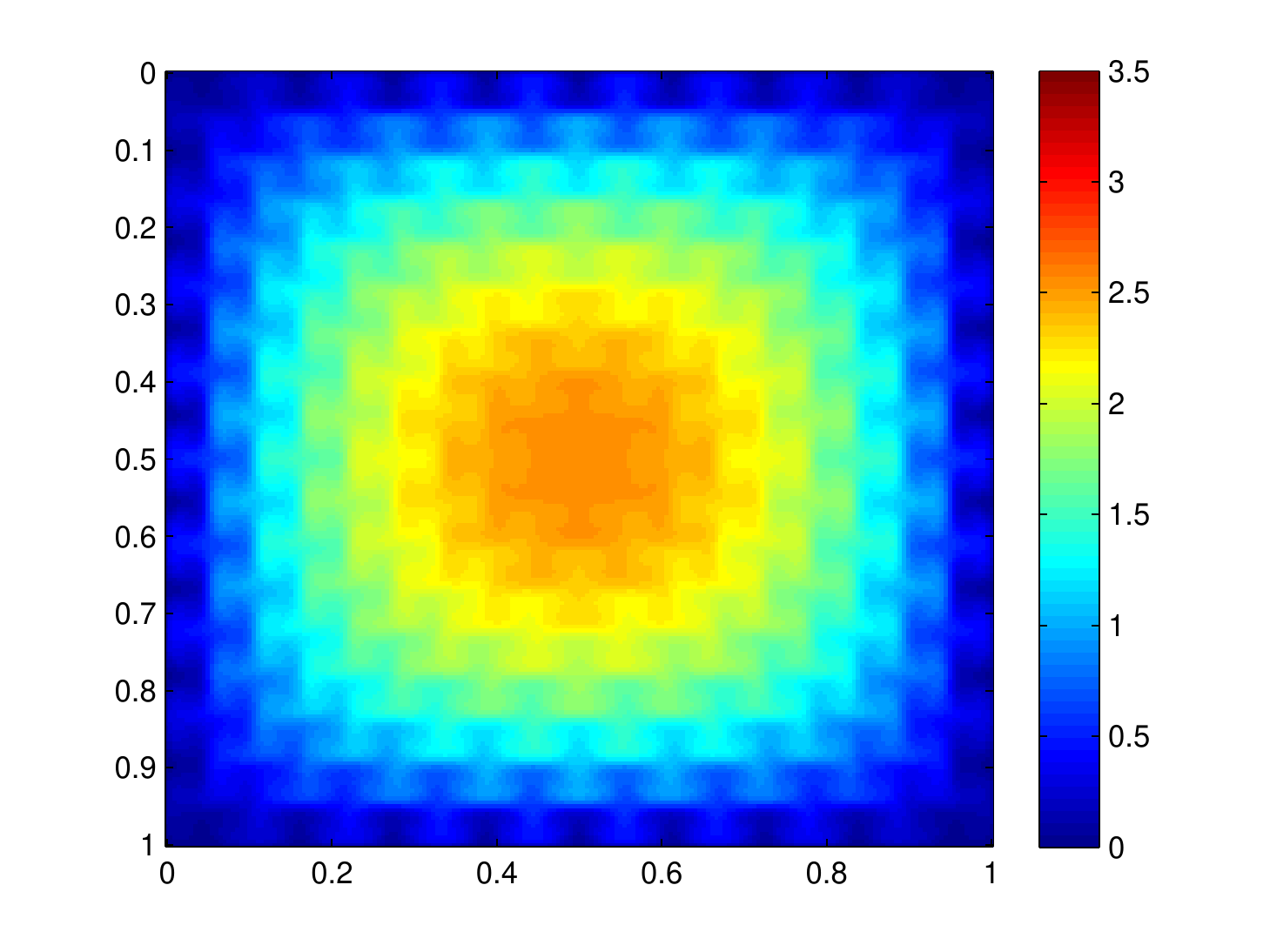}
  }
  \caption{Reference solutions for Example 1.}
  \label{fig:solns_ex1}
\end{figure}

\begin{figure}[!ht]
  \centering
  \subfigure[ $\alpha=2$]{\label{fig:solns_ex2L}
    \includegraphics[width = 0.475\textwidth]{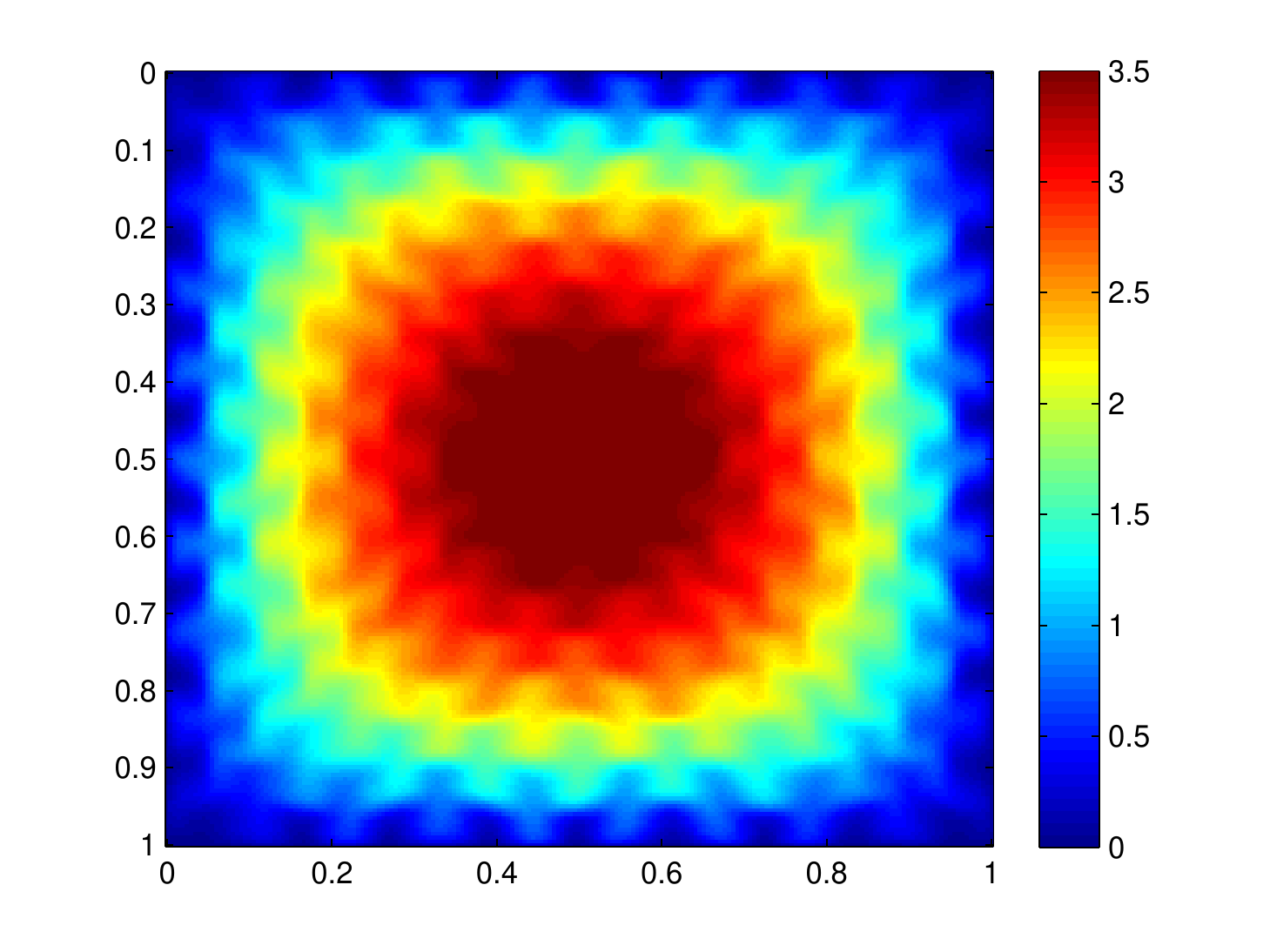}
   }
  \subfigure[ $\alpha=4$]{\label{fig:solns_ex2R}
     \includegraphics[width = 0.475\textwidth]{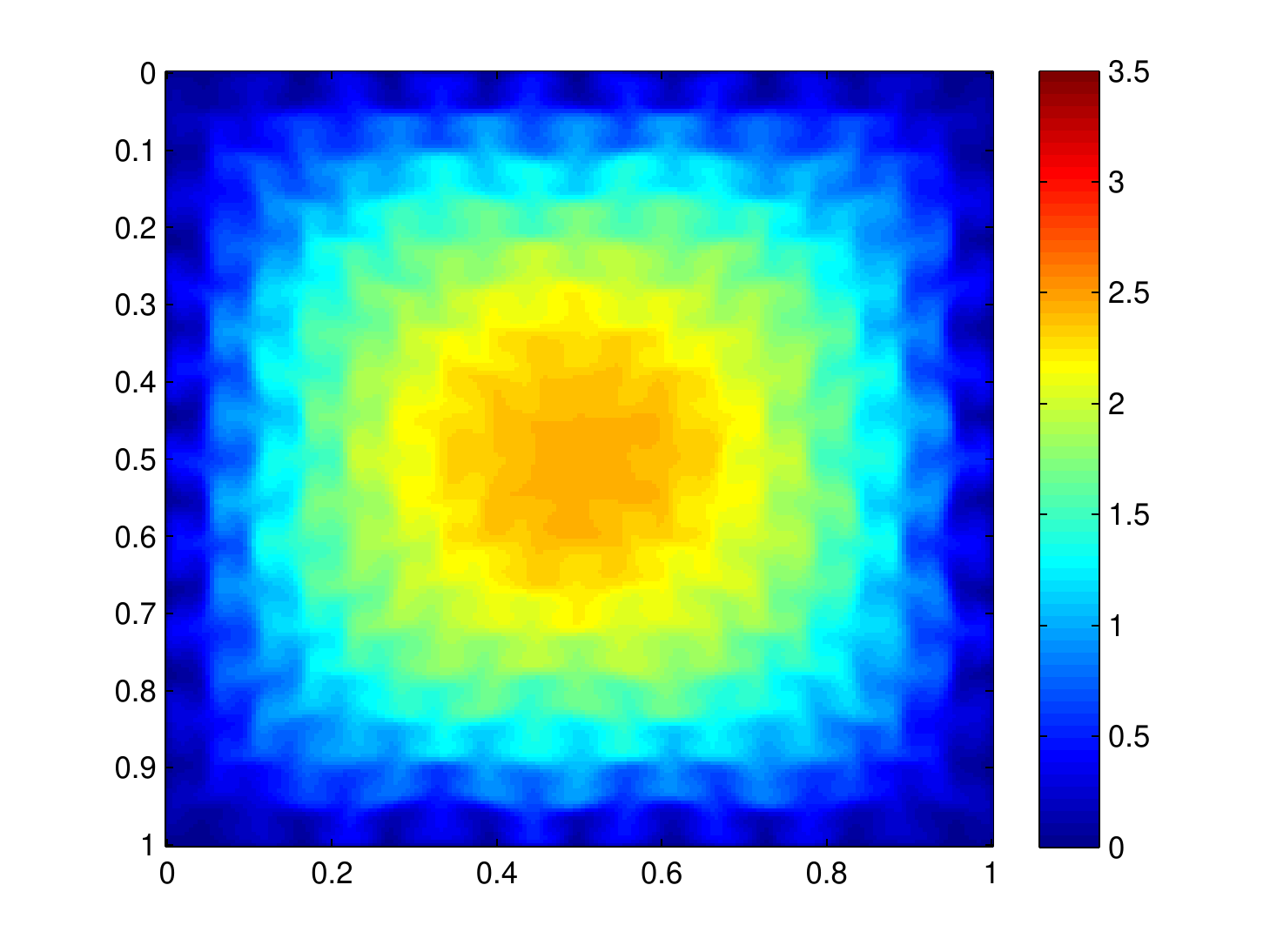}
  }
  \caption{Reference solutions for Example 2.}
  \label{fig:solns_ex2}
\end{figure}

In this section, we present representative numerical examples. In all
our examples, $\{\chi_{i}\}$ is a multiscale partition of unity. In
each coarse space, we compare the $l_2$ projection error and the $l_2$
error for the multiscale solution. For simplicity, we refer to ``the
multiscale error'' as the error between the multiscale solution and
the exact solution, and ``the projection error'' as the error between
the exact solution and its $l_2$ projection onto the span of the
coarse trial space.  We also assume $\kappa$ is a constant and $b$ is
a multiscale field. In particular, the velocity fields contain
oscillations and cells (eddies, separatrices and/or layers) within a
single coarse block of the discretization and, thus, we do not have a
single streamline direction per coarse block.  Fully resolved velocity
solutions are shown in Figures~\ref{fig:solns_ex1}
to~\ref{fig:solns_ex3} and in Figure~\ref{fig:solns_ex5}. The method
can easily handle multiscale diffusion coefficients.  The fine-grid
problem is always chosen such that the local Peclet number is about
$1$ ensuring a stable fine discretization. All coarse discretizations
have a Peclet number at least an order of magnitude larger than
$1$. We analyze the performance of the trial and test
spaces proposed in the previous section. We pay special attention to
the effect of eigenvalue problem on the performance of the discrete
system and discuss this for each example.

\subsection*{Example 1}

These first numerical examples are defined by the following diffusion and
convection coefficients and right-hand side,
\begin{align*}
  \kappa & =\dfrac{1}{100},\\
  b&=\alpha \left(\begin{array}{c}
      +\sin(18\pi x)\cos(18\pi y)\\
      -\cos(18\pi x)\sin(18\pi y)
\end{array}\right),\\
f & =1.
\end{align*}
\begin{table}[!ht]
\centering
\begin{tabular}{ccccc}
\toprule
\#basis & \multicolumn{4}{c}{$L_{2}$ (projection error)}\tabularnewline
 (trial, test)
& \multicolumn{2}{c}{Eigenproblem 1}& \multicolumn{2}{c}{Eigenproblem 2}\tabularnewline
\midrule
&$\alpha=2$&$\alpha=4$&$\alpha=2$&$\alpha=4$\tabularnewline
\midrule
(1,1) & 8.56\%& 7.06\%& 11.94\%& 9.60\%\tabularnewline
(1,3) & 3.22\%& 4.96\%& 4.74\%& 4.48\%\tabularnewline
(1,5) & 2.85\%& 4.74\%& 2.90\%& 5.02\%\tabularnewline
(1,7) & 2.85\%(2.85\%)& 3.64\%(3.52\%)& 2.85\%(2.85\%)& 3.55\%(3.52\%)\tabularnewline
\midrule
(3,1) & 9.00\%& 7.58\%&11.95\%& 8.86\%\tabularnewline
(3,3) & 3.12\%& 5.22\%& 5.01\%& 3.96\%\tabularnewline
(3,5) & 2.61\% & 3.96\%& 2.70\% & 4.83\%\tabularnewline
(3,7) & 2.60\%(2.60\%)& 3.41\%(3.21\%)& 2.61\%(2.60\%)& 3.25\%(3.21\%)\tabularnewline
\midrule
(5,1) & 8.65\%& 7.88\%& 12.80\%& 9.08\%\tabularnewline
(5,3) & 2.72\%& 4.97\%& 4.69\%& 3.35\%\tabularnewline
(5,5) & 2.31\%& 3.62\%& 2.37\%& 3.99\%\tabularnewline
(5,7) & 2.31\%(2.31\%)& 2.89\%(2.77\%)& 2.31\%(2.31\%)& 2.79\%(2.77\%)\tabularnewline
\bottomrule
\end{tabular}
\caption{Errors for test space derived using Eigenproblems 1 and 2 for
Example 1. Coarse and fine mesh sizes are $H=1/10$ and $h=1/200$,
respectively. The projection errors are in parentheses.}
\label{table:ex1}
\end{table}

This velocity field has a cellular structure with several eddies and
separatices. In the simulations $\alpha$ takes values of 2 and
4. Figure~\ref{fig:solns_ex1} depicts well-resolved fine-scale
solutions for the chosen values of $\alpha$. In both cases, we take
the coarse mesh size to be $H=1/10$, while the fine mesh size to be
$h=1/200$.

\begin{table}[!ht]
\centering
\begin{tabular}{ccc}
\toprule
\#basis test & \multicolumn{2}{c}{$\min\{\lambda_{L_{i}+1}\}$}\tabularnewline
\midrule
&$\alpha=2$&$\alpha=4$\tabularnewline
\midrule
1 & 0.3445 & 0.3693\tabularnewline
3 & 0.7273 & 0.7707\tabularnewline
5 & 0.9542 & 0.9514\tabularnewline
7 & 0.9908 & 0.9919\tabularnewline
\bottomrule
\end{tabular}
\caption{Minimum eigenvalue for the test space constructed using
   the Eigenproblem 2 (minimal energy test functions) for Example 1. Coarse and fine
  mesh  sizes are $H=1/10$ and $h=1/200$,  respectively.}
\label{tab:eigex1}
\end{table}

Table~\ref{table:ex1} shows the impact of increasing the number of
coarse basis functions as well as how the system converges as we
increase the number of test functions included per coarse block
edge. The table shows the evolution of the multiscale error as we
increase the number of test functions for different numbers of trial
functions in each coarse block. Each column is labeled by its
corresponding value of $\alpha$.  Table~\ref{table:ex1} shows the
performance of the reduced-dimensional test space constructed using
the first and second eigenvalue problems we describe in
Section~\ref{sec:TestEigenproblem}. This table shows that $7$ test
functions per edge of the coarse mesh are enough to deliver similar
multiscale and projection errors irrespective of $\alpha$ (i.e.,
coarse scale Peclet number) and the number of coarse basis functions
used in each coarse block. In fact, these errors are similar even when
the number of test functions is $5$. Table~\ref{tab:eigex1} shows the
evolution of the minimum eigenvalue for the test space constructed
using the Eigenproblem 2 (minimal energy test functions) of
Section~\ref{sec:TestEigenproblem}.  As it follows from the theory,
for a rich enough test space with a sufficient number of multiscale
test functions, multiscale and projection errors converge. The
eigenvalue behavior shown in Table~\ref{tab:eigex1} and the
convergence shown in Table~\ref{table:ex1} verifies that when the
minimum eigenvalue is close to $1$, the multiscale error converges to
the projection error.

\subsection*{Example 2}

 We consider the following right-hand side, and diffusion and
convection coefficients:
\begin{align*}
  \kappa  &=1/100,\\
  b&=\alpha\left(\begin{array}{c}
                   +\sin(18\pi x) \cos(18\pi y) + \delta \cos(18 \sqrt{2} \pi x) \sin(18 \sqrt{2} \pi y)\\
                   -\cos(18\pi x)\sin(18\pi y -\delta \sin(18\sqrt{2} \pi x)\sin(18\sqrt{2} \pi y)
                 \end{array}\right),\\
  f & =1.
\end{align*}
This velocity field has a cellular structure with eddies and channels,
as can be seen in the fine-scale solutions depicted in
Figure~\ref{fig:solns_ex2}. These channels introduce global
effects. $\alpha$ takes values of 2 and 4, as in our previous
example. $\delta$ takes value of $\sqrt{2}/4$.   For both
cases, we take the coarse mesh size to be $H=1/10$, while the fine
mesh size to be $h=1/200$.

\begin{table}[!ht]
\centering
\begin{tabular}{ccccc}
\toprule
\#basis & \multicolumn{4}{c}{$L_{2}$ (projection error)}\tabularnewline
 (trial, test)
& \multicolumn{2}{c}{Eigenproblem 1}& \multicolumn{2}{c}{Eigenproblem 2}\tabularnewline
\midrule
&$\alpha=2$&$\alpha=4$&$\alpha=2$&$\alpha=4$\tabularnewline
\midrule
(1,1) & 9.49\%& 8.57\%& 12.27\% & 9.99\%\tabularnewline
(1,3) & 3.11\%& 5.04\%& 4.62\%& 4.22\%\tabularnewline
(1,5) & 2.90\%& 4.73\%& 2.95\% & 4.62\%\tabularnewline
(1,7) & 2.90\%(2.90\%)& 3.87\%(3.67\%) & 2.90\%(2.90\%) & 3.70\%(3.67\%)\tabularnewline
\midrule
(3,1) & 10.04\%& 9.09\%& 12.35\% & 9.33\%\tabularnewline
(3,3) & 2.97\%&4.50\%& 5.09\%& 3.73\%\tabularnewline
(3,5) & 2.65\% & 3.85\%& 2.75\% & 4.52\%\tabularnewline
(3,7) & 2.65\%(2.64\%)& 3.56\%(3.33\%) & 2.65\%(2.64\%) & 3.38\%(3.33\%)\tabularnewline
\midrule
(5,1) & 9.22\%&9.25\%& 13.37\% & 9.98\%\tabularnewline
(5,3) & 2.67\%& 4.18\%& 4.60\%& 3.26\%\tabularnewline
(5,5) & 2.36\%& 3.43\%& 2.42\%& 3.84\%\tabularnewline
(5,7) & 2.36\%(2.36\%)& 3.18\%(2.84\%) & 2.36\%(2.36\%) & 2.88\%(2.84\%)\tabularnewline
\bottomrule
\end{tabular}
\caption{Errors for test space derived using Eigenproblems 1 and 2 for
Example 2, with $\delta=\sqrt{2}/4$. Coarse and fine mesh sizes
are $H=1/10$ and $h=1/200$, respectively. The projection errors are in
parentheses.}
\label{table:ex2}
\end{table}

\begin{table}[!ht]
\centering
\begin{tabular}{ccc}
\toprule
\#basis test & \multicolumn{2}{c}{$\min\{\lambda_{L_{i}+1}\}$}\tabularnewline
\midrule
&$\alpha=2$&$\alpha=4$\tabularnewline
\midrule
1 & 0.3551 & 0.2967\tabularnewline
3 & 0.7289 & 0.6761\tabularnewline
5 & 0.9510 & 0.9331\tabularnewline
7 & 0.9909 & 0.9813\tabularnewline
\bottomrule
\end{tabular}
\caption{Minimum eigenvalue for the test space constructed using
   the Eigenproblem 2 (minimal energy test functions) for Example 2
    for $\delta=\sqrt{2}/4$. Coarse and fine
  mesh  sizes are $H=1/10$ and $h=1/200$,  respectively.}
\label{tab:eigex2}
\end{table}

In Table~\ref{table:ex2}, we increase the number of test functions and
consider different numbers of trial functions for both eigenvalue
problems described in Section~\ref{sec:TestEigenproblem}. As
before, we observe that for $7$ test functions, the
multiscale and projection errors converge, while for $5$ test
functions per edge, the errors are close.  As we increase, the dimension of the
coarse trial space, we observe a similar behavior. The eigenvalue
behavior shows (see Table~\ref{tab:eigex2}) that when the eigenvalue
is close to $1$, the multiscale and projection errors converge.

\begin{figure}[!ht]
  \centering
  \subfigure[ $\alpha=1/1000$]{\label{fig:solns_ex3L}
    \includegraphics[width = 0.475\textwidth]{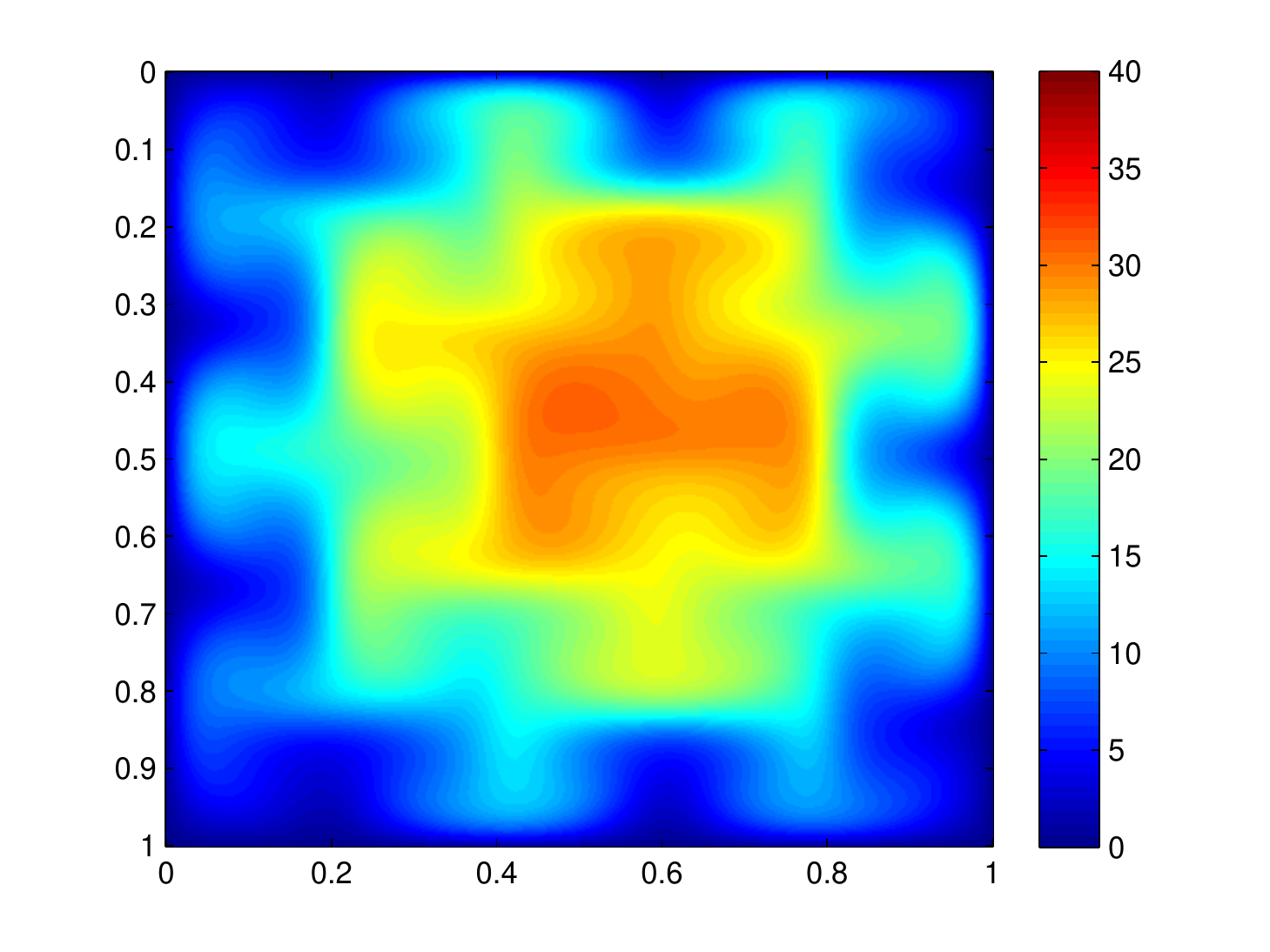}
   }
  \subfigure[ $\alpha=1/2000$]{\label{fig:solns_ex3R}
     \includegraphics[width = 0.475\textwidth]{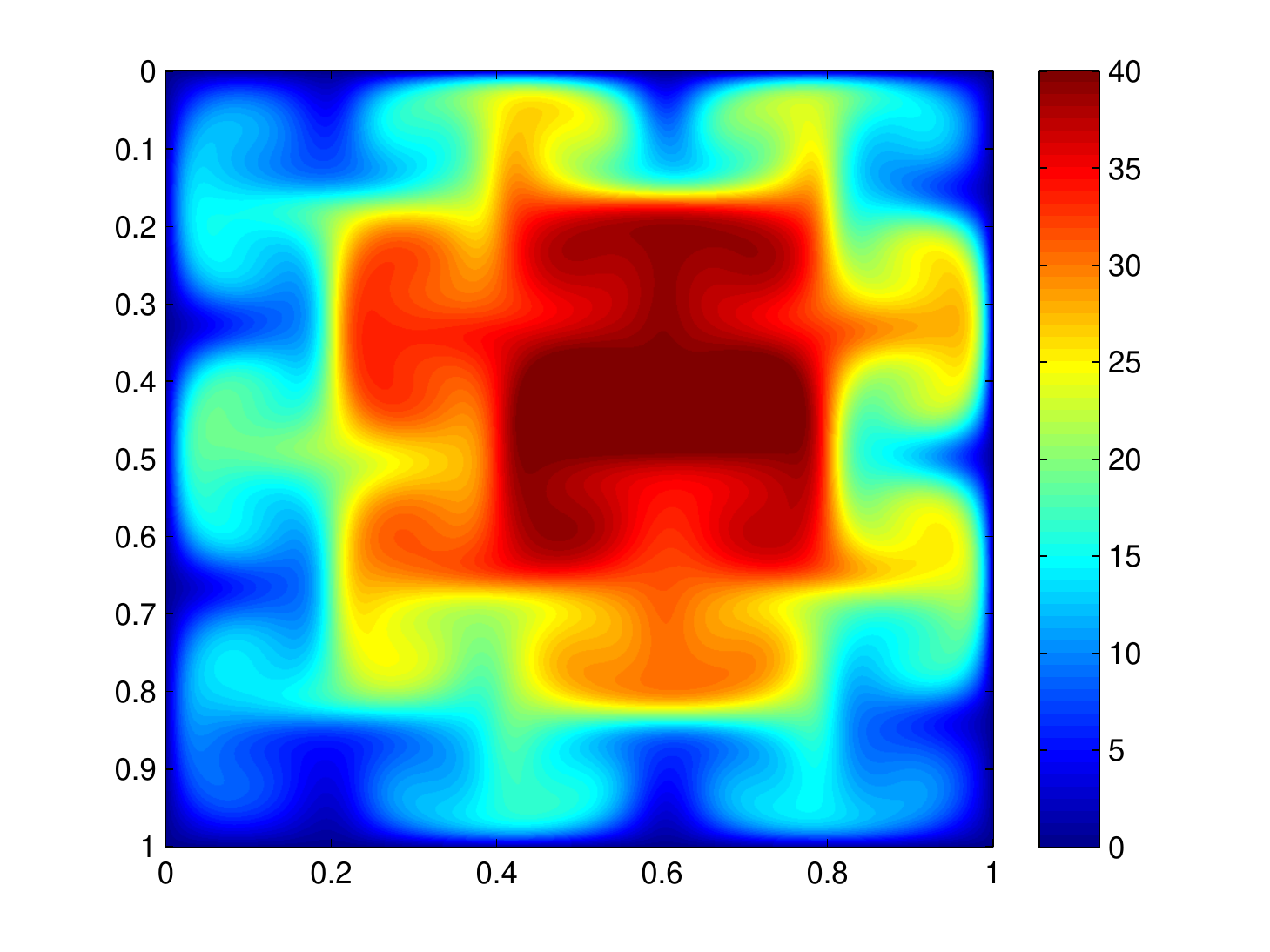}
  }
  \caption{Reference solutions for Example 3.}
  \label{fig:solns_ex3}
\end{figure}

\subsection*{Example 3}

\begin{table}[!ht]
\centering
\begin{tabular}{ccccc}
\toprule
\#basis & \multicolumn{4}{c}{$L_{2}$ (projection error)}\tabularnewline
 (trial, test)
& \multicolumn{2}{c}{Eigenproblem 1}& \multicolumn{2}{c}{Eigenproblem 2}\tabularnewline
\midrule
&$\alpha=1/1000$&$\alpha=1/2000$&$\alpha=1/1000$&$\alpha=1/2000$\tabularnewline
\midrule
(1,1) & 9.59\%& 32.36\%& 19.09\% &34.55\%\tabularnewline
(1,3) & 7.84\%& 13.70\% & 15.58\%&  32.00\%\tabularnewline
(1,5) & 7.83\%& 12.10\% & 9.10\% & 22.53\%\tabularnewline
(1,7) & 7.83\%(7.83\%)& 11.85\%(11.48\%)& 7.95\%(7.83\%)  & 14.70\%(11.48\%)\tabularnewline
\midrule
(3,1) &  6.01\%& 16.24\%& 16.70\%& 37.12\%\tabularnewline
(3,3) & 4.89\% &  9.07\%& 7.93\%& 18.82\%\tabularnewline
(3,5) & 4.88\% &  7.98\%& 5.33\% &11.90\%\tabularnewline
(3,7) &  4.88\%(4.88\%)& 8.06\%(7.58\%) & 4.98\%(4.88\%) & 9.71\%(7.58\%)\tabularnewline
\midrule
(5,1) & 4.70\%&14.77\%&  14.00\%&  31.46\%\tabularnewline
(5,3) & 3.74\%& 6.61\%& 4.79\%& 13.66\%\tabularnewline
(5,5) & 3.72\%& 6.23\%& 3.80\%& 7.55\%\tabularnewline
(5,7) & 3.72\%(3.72\%)& 6.21\%(6.15\%) & 3.74\%(3.72\%) & 6.53\%(6.15\%)\tabularnewline
\bottomrule
\end{tabular}
\caption{Errors for test space derived using Eigenproblems 1 and 2 for
  Example 3. Discrete parameters used are: $\alpha=1/1000$, $H=1/10$,
  $h=1/400$ and $\alpha=1/2000$, $H=1/10$, $h=1/800$. The
  projection errors are in parentheses.}
\label{table:ex3}
\end{table}

\begin{table}[!ht]
\centering
\begin{tabular}{ccc}
\toprule
\#basis test & \multicolumn{2}{c}{$\min\{\lambda_{L_{i}+1}\}$}\tabularnewline
\midrule
&$\alpha= 1/1000$&$\alpha= 1/2000$\tabularnewline
\midrule
1 & 0.3547 & 0.3068\tabularnewline
3 & 0.7497 & 0.6304\tabularnewline
5 & 0.9546 & 0.8718\tabularnewline
7 & 0.9952 & 0.9754\tabularnewline
\bottomrule
\end{tabular}
\caption{Minimum eigenvalue for test space derived using Eigenproblem
  2 for Example 3. Discrete parameters used are: $\alpha=1/1000$, $H=1/10$,
  $h=1/400$ and $\alpha=1/2000$, $H=1/10$, $h=1/800$.  }
\label{tab:eigex3}
\end{table}

We consider the following diffusion and convection coefficients and
right-hand side.
\begin{align*}
  \kappa & =\alpha,\\
  b&= \left(\begin{array}{c}
                    -\frac{\partial H}{\partial y}\vspace{.05in}\\
                    +\frac{\partial H}{\partial x}
\end{array}\right),\\
  f & =1,
\end{align*}
where \[H=(\sin(5\pi x)\sin(6\pi y)/(60\pi))+0.005(x+y).\]
This velocity field again has a cellular structure with eddies and
channels, as the fine-scale solutions show in
Figure~\ref{fig:solns_ex3}. In this example, $\alpha$ is a diffusion
coefficient and we take $\alpha=1/1000$ and $\alpha=1/2000$.  For both
cases, we take the coarse mesh size to be $H=1/10$, while the fine
mesh sizes are set to $h=1/400$ and $h=1/800$ for $\alpha=1/1000$ and
$\alpha=1/2000$, respectively.

Tables~\ref{table:ex3} and~\ref{tab:eigex3} show a similar behavior to
the one discussed in the previous two examples. That is,
table~\ref{table:ex3} shows that if only one test function is chosen
for $\alpha=1/2000$, the error is about $35$\% (when the number of
test functions coarse edge is $1$ and the number of trial functions
per coarse block is $1$).  These errors rapidly drop to about the
projection error as we increase the dimension of the test
space. Similar behavior is observed for both eigen-constructions of
the test space as we refine the coarse trial space.

\begin{table}[!ht]
\centering

\begin{tabular}{cc}
\toprule
\#basis  (trial, test)& \multicolumn{1}{c}{$L_{2}$ (projection
  error)}\tabularnewline
\midrule
(1,1) & 20.12\%\tabularnewline
(1,3) & 11.25\%\tabularnewline
(1,5) & 4.03\%\tabularnewline
(1,7) & 3.93\%(3.93\%)\tabularnewline
\midrule
(3,1) & 19.99\%\tabularnewline
(3,3) & 13.02\%\tabularnewline
(3,5) & 3.31\%\tabularnewline
(3,7) & 3.23\%(3.23\%)\tabularnewline
\midrule
(5,1) & 13.93\%\tabularnewline
(5,3) & 9.14\%\tabularnewline
(5,5) & 2.90\%\tabularnewline
(5,7) & 2.74\%(2.70\%)\tabularnewline
\bottomrule
\end{tabular}
\caption{Errors for test space derived using Eigenproblem 2 for
  Example 4. Coarse and fine mesh  sizes are $H=1/10$ and $h=1/200$,
  respectively.  The projection errors are in parentheses.}
\label{table:ex4}
\end{table}

\subsection*{Example 4}
As we remove the eddies from the flow field and make the flow more
channelized, the multiscale error grows.  To expose this behavior, we
take the velocity field to be
\begin{align*}
  \kappa & =1,\\
  b&=200 \left( \begin{array}{c}
                  \sin(18\sqrt{2} \pi y)\\
                  0 \end{array}\right ),\\
  f & =1,
\end{align*}
which corresponds to solving the flow equations with a channelized
permeability field.  The numerical results are presented in
Table~\ref{table:ex4} (the mesh sizes for the coarse and fine spaces
are $H=1/10$ and $h=1/200$).  We observe that the multiscale error is
$20.12$\% for one trial function per coarse block and one test
function per coarse interface, while the error reduces to
the projection error of $3.93$\% when we select $7$ test functions
interface. As before, for $5$ trial functions per coarse block, it
takes $7$ test functions per coarse interface to reduce the error to
the projection error from $13.93$\%.  For this discrete problem setup,
the smallest eigenvalue is $0.9952$ for $7$ test functions per edge when
minimal energy functions are used (Eigenproblem 2 in
Section~\ref{sec:TestEigenproblem}).

\begin{figure}[!ht]
\centering
\includegraphics[scale=0.65]{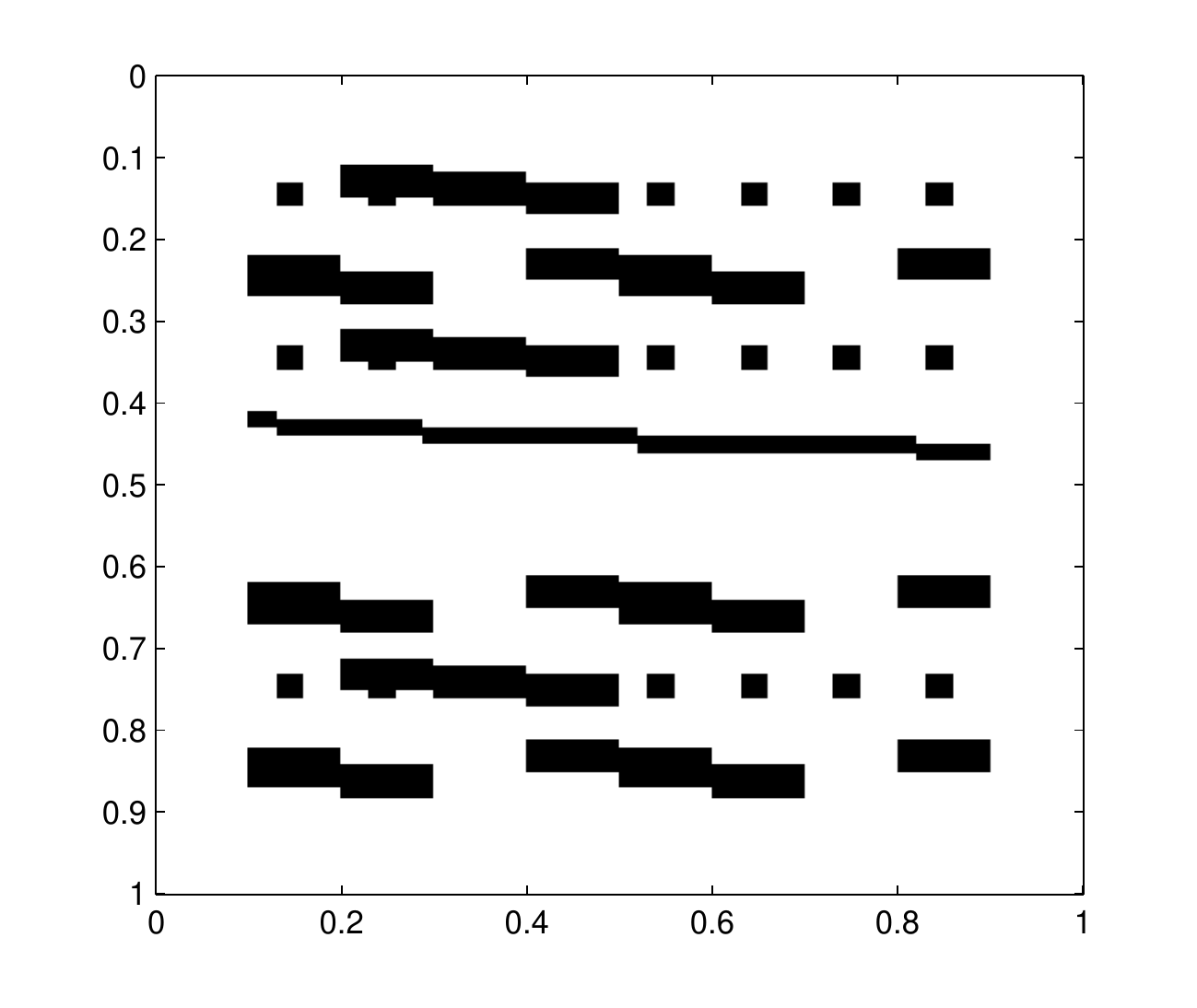}
\caption{Permeability field used to compute the transport velocity
  field in Example 5. The black region corresponds to the permeability $500$ and the white region corresponds to the permeability $1$.}
\label{fig:kappa}
\end{figure}

\begin{figure}[!ht]
  \centering
  \subfigure[ $\alpha=1/250$]{\label{fig:solns_ex5L}
    \includegraphics[width = 0.475\textwidth]{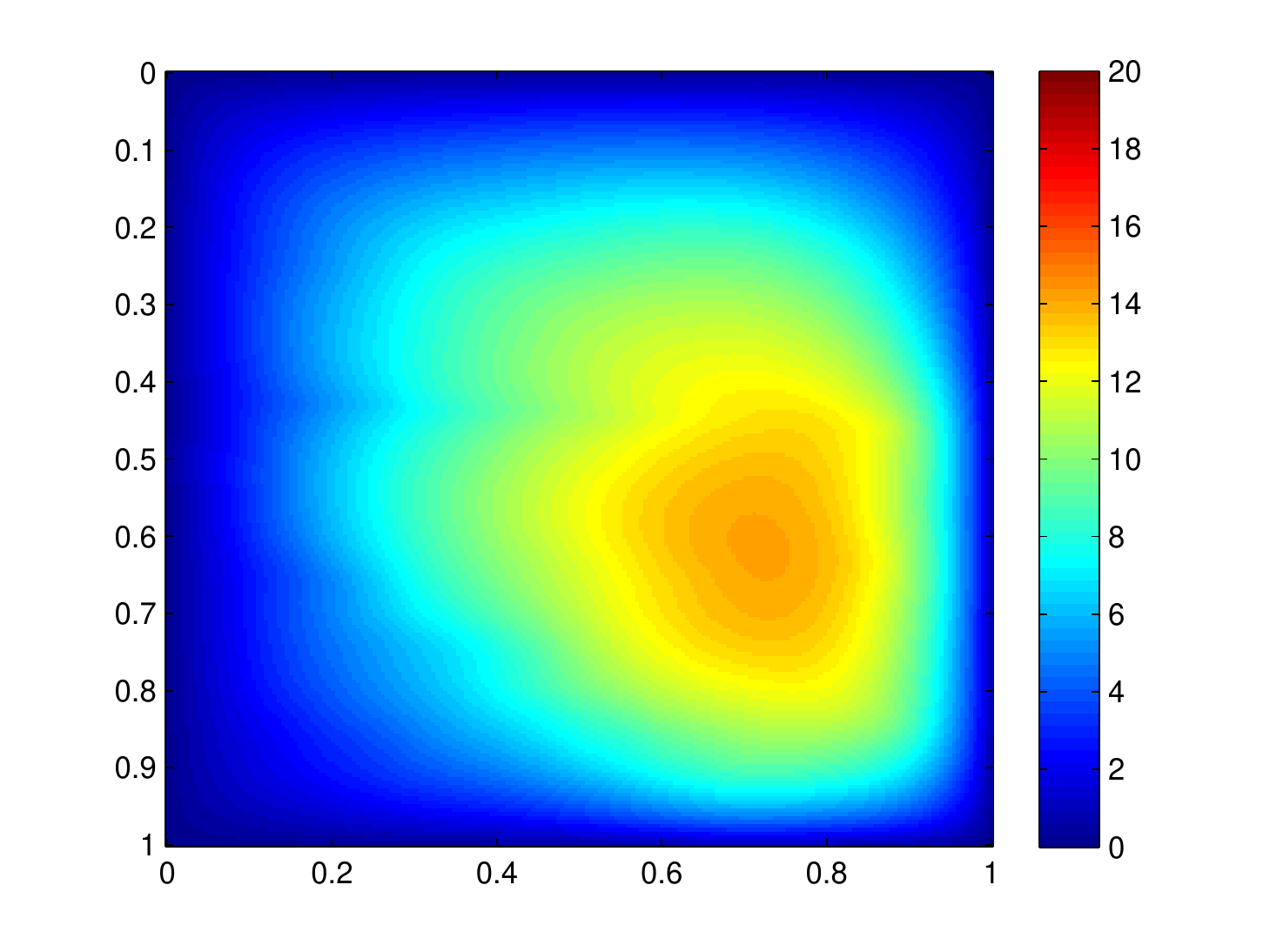}
   }
  \subfigure[ $\alpha=1/500$]{\label{fig:solns_ex5R}
     \includegraphics[width = 0.475\textwidth]{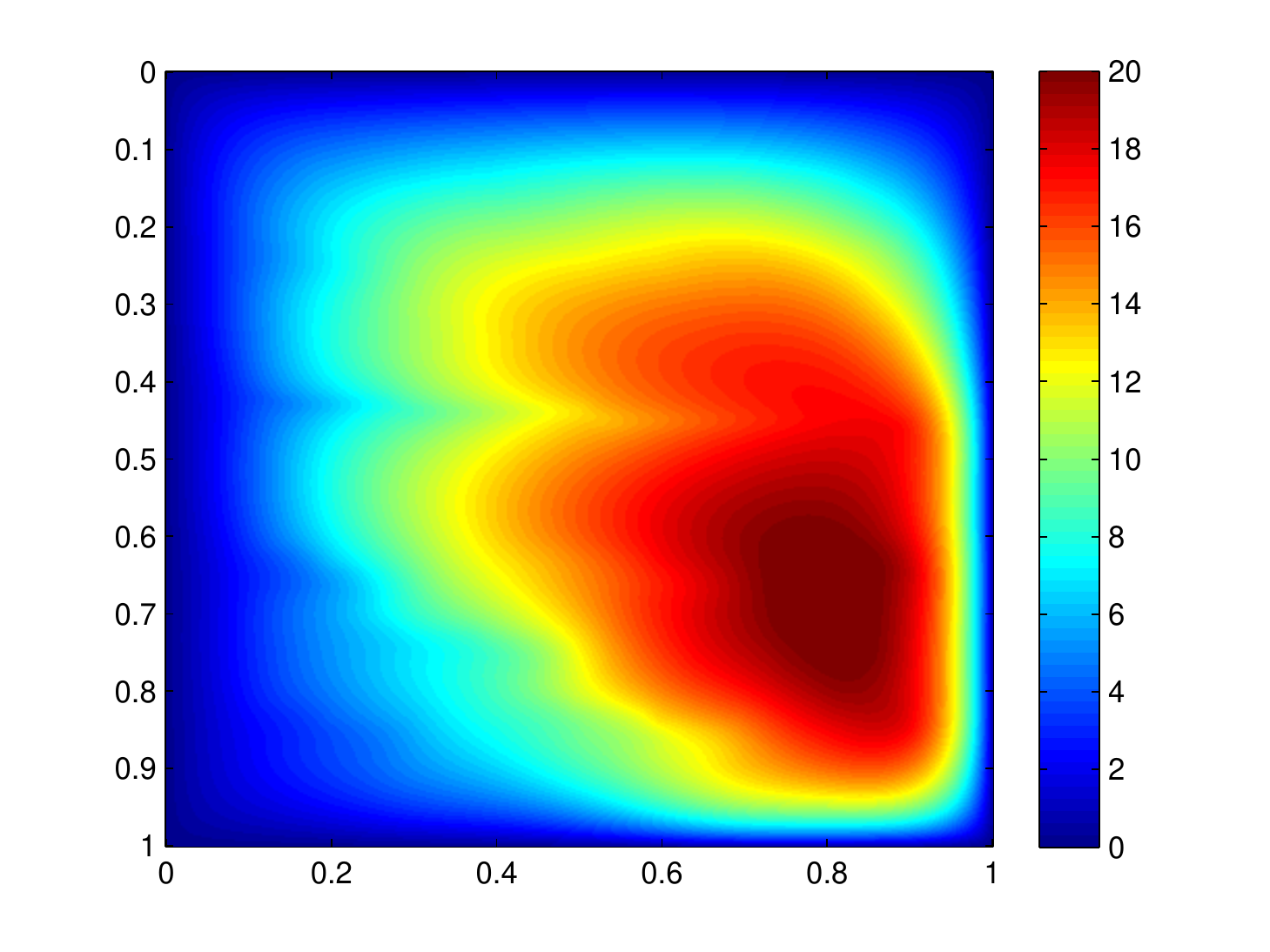}
  }
  \caption{Reference solutions for Example 5.}
  \label{fig:solns_ex5}
\end{figure}

\subsection*{Example 5}

\begin{table}[!ht]
\centering
\begin{tabular}{ccccc}
\toprule
\#basis & \multicolumn{4}{c}{$L_{2}$ (projection error)}\tabularnewline
 (trial, test)
& \multicolumn{2}{c}{Eigenproblem 1}& \multicolumn{2}{c}{Eigenproblem 2}\tabularnewline
\midrule
&$\alpha=1/250$&$\alpha=1/500$&$\alpha=1/250$&$\alpha=1/500$\tabularnewline
\midrule
(1,1) & 2.11\%& 4.64\%& 2.24\%& 4.56\%\tabularnewline
(1,3) & 2.08\%& 4.59\% & 2.09\%&  4.24\%\tabularnewline
(1,5) & 2.07\%& 4.45\% & 2.07\%& 4.15\%\tabularnewline
(1,7) & 2.07\%(2.07\%)& 4.23\%(4.04\%)&2.07\%(2.07\%)  & 4.19\%(4.04\%)\tabularnewline
\midrule
(3,1) & 1.01\%& 1.98\% & 1.68\%& 3.57\%\tabularnewline
(3,3) & 0.99\% &  1.93\%& 1.07\%& 2.36\%\tabularnewline
(3,5) & 0.99\% &  1.93\%& 1.00\%& 2.05\%\tabularnewline
(3,7) & 0.99\%(0.99\%)& 1.93\%(1.91\%) & 0.99\%(0.99\%)& 2.01\%(1.91\%)\tabularnewline
\midrule
(5,1) & 0.85\%&1.70\%&  1.64\%&  4.12\%\tabularnewline
(5,3) & 0.75\%& 1.44\%& 0.84\%& 1.91\%\tabularnewline
(5,5) & 0.75\%& 1.44\%& 0.76\%& 1.53\%\tabularnewline
(5,7) & 0.75\%(0.75\%)& 1.44\%(1.42\%) & 0.75\%(0.75\%)& 1.49\%(1.42\%)\tabularnewline
\bottomrule
\end{tabular}
\caption{Errors for test space derived using Eigenproblems 1 amd 2 for
Example 5. Discrete parameters used are: $\alpha=1/250$, $H=1/10$,
$h=1/200$ and $\alpha=1/500$, $H=1/10$, $h=1/400$. The projection
errors are in parentheses.}
\label{table:ex5}
\end{table}

We consider the following diffusion and convection coefficients, and
right-hand side.
\begin{align*}
  \kappa & =\alpha,\\b&=\kappa\nabla p\\
  f & =1,
\end{align*}
where the velocity field solves this flow equation
\begin{align*}
  -\nabla\cdot(\kappa\nabla p) & =0\\
  p|_{\partial\Omega} & =xy
\end{align*}
Figure~\ref{fig:kappa} shows the permeability field used in the above
equation. The resulting velocity field contains channels with variable
velocity in each coarse region. Figure~\ref{fig:solns_ex5} shows the
fine-scale structure of the fully resolved velocity field.  In this
case $\alpha$ is a diffusion coefficient and takes values $1/250$ and
$1/500$.  In both cases, we take the coarse mesh size is set to
$H=1/10$, while the fine mesh sizes are $h=1/200$ and $h=1/400$ for
$\alpha=1/250$ and $\alpha=1/500$, respectively.

Tables~\ref{table:ex5} and~\ref{tab:eigex5} show a similar
behavior to that observed in the prior examples. That is, the
multiscale error converges to the projection error as we increase the
number of test functions per coarse edge for either eigenvalue problem
and for any number of coarse functions in each coarse block.

\begin{table}[!ht]
\centering
\begin{tabular}{ccc}
\toprule
\#basis test & \multicolumn{2}{c}{$\min\{\lambda_{L_{i}+1}\}$}\tabularnewline
\midrule
&$\alpha= 1/1000$&$\alpha= 1/2000$\tabularnewline
\midrule
1 & 0.4106 & 0.3544\tabularnewline
3 & 0.8592 & 0.7583\tabularnewline
5 & 0.9828 & 0.9535\tabularnewline
7 & 0.9985 & 0.9919\tabularnewline
\bottomrule
\end{tabular}
\caption{Minimum eigenvalue for test space derived using Eigenproblem
2 for Example 3. Discrete parameters used
are:$\alpha=1/250,H=1/10,h=1/200$ and $\alpha=1/500,H=1/10,h=1/400$}
\label{tab:eigex5}
\end{table}

\begin{table}[!ht]
\centering
\begin{tabular}{cccccc}
  \toprule
%&& \multicolumn{4}{c}{$L_{2}$ (projection error)}\tabularnewline
   \#basis &\#iter & \multicolumn{2}{c}{Eigenproblem 1}& \multicolumn{2}{c}{Eigenproblem 2}\tabularnewline
   (trial, test) &&$\alpha=2$&$\alpha=4$&$\alpha=2$&$\alpha=4$\tabularnewline
  \midrule
  & 0 & 8.56\%&7.06\% & 11.94\% & 9.60\%\tabularnewline
  (1,1)&1 & 2.89\%& 3.79\%& 2.96\% & 4.30\%\tabularnewline
  &2 & 2.85\%(2.85\%) & 3.52\%(3.52\%) & 2.85\%(2.85\%) & 3.52\%(3.52\%)\tabularnewline
  \midrule
  &0 & 3.22\%& 4.96\%& 4.74\% & 4.48\%\tabularnewline
  (1,3)&1 & 2.85\%& 3.54\%& 2.86\% & 3.63\%\tabularnewline
  &2 & 2.85\%(2.85\%) & 3.52\%(3.52\%) & 2.85\%(2.85\%) & 3.52\%(3.52\%)\tabularnewline
  \midrule
  &0 & 8.65\%& 7.88\% & 12.80\% & 9.08\%\tabularnewline
  (5,1)& 1 & 2.33\%& 2.97\% & 2.58\% & 3.29\%\tabularnewline
  &2 & 2.31\%(2.31\%) & 2.77\%(2.77\%) & 2.32\%(2.31\%) & 2.78\%(2.77\%)\tabularnewline
  \midrule
  &0 & 2.72\%& 4.97\% & 4.69\% & 3.35\%\tabularnewline
  (5,3)& 1 & 2.31\%& 2.79\% & 2.33\% & 2.81\%\tabularnewline
  & 2 & 2.31\%(2.31\%) & 2.77\%(2.77\%) & 2.31\%(2.31\%) & 2.77\%(2.77\%)\tabularnewline
  \bottomrule
\end{tabular}
\caption{Error evolution as online basis functions are added to the
system (test space derived using Eigenproblems 1 and 2 for
Example~1). Coarse and fine mesh sizes are $H=1/10$ and $h=1/200$,
respectively.  The projection errors are in parentheses.}
\label{tab:onlineex1}
\end{table}

\subsection{Numerical Result for online test basis enrichment}

In this section, we present some numerical results, which use online
test basis functions to stabilize the system. In
Table~\ref{tab:onlineex1}, we show the convergence history for the
online test basis enrichment for Example 1, while in
Table~\ref{tab:onlineex4}, we show the convergence history of the
online test basis enrichment for Example 4. In these two cases, with only
one iteration, the multiscale error becomes similar to the projection
error. In the second iteration, the multiscale error converges to
the projection error.

\begin{table}[!ht]
\centering
\begin{tabular}{cccccc}
  \toprule
%&& \multicolumn{4}{c}{$L_{2}$ (projection error)}\tabularnewline
   \#basis &\#iter & \tabularnewline
   (trial, test) &  \tabularnewline
  \midrule
  & 0 & 20.12\%\tabularnewline
  (1,1)&1 & 3.93\%\tabularnewline
  &2 & 3.92\%(3.92\%) \tabularnewline
  \midrule
  &0 & 11.15\%\tabularnewline
  (1,3)&1 & 3.92\%\tabularnewline
  &2 & 3.92\%(3.92\%) \tabularnewline
  \midrule
  &0 & 13.93\%\tabularnewline
  (5,1)& 1 &3.24\%\tabularnewline
  &2 & 2.72\%(2.70\%) \tabularnewline
  \midrule
  &0 & 9.14\%\tabularnewline
  (5,3)& 1 & 2.74\%\tabularnewline
  & 2 & 2.70\%(2.70\%) \tabularnewline
  \bottomrule
\end{tabular}
\caption{Error evolution as online basis functions are added to the
system (test space derived using Eigenproblem 2 for
Example~4). Coarse and fine mesh sizes are $H=1/10$ and $h=1/200$,
respectively.  The projection errors are in parentheses.}
\label{tab:onlineex4}
\end{table}

\section{Conclusions}

In this paper, we study multiscale methods for convection-dominated
diffusion with heterogeneous convective velocity fields. This
stabilization generalizes the approaches described in~\cite{
  demkowicz2014overview} to multiscale problems. To construct this
stabilization we reformulate overall problem in mixed form.  The
auxiliary variable we introduce plays the role of the test
function. We describe the multiscale spaces we use for the test and
trial spaces, which are built using GMsFEM framework. First, we
construct snapshots spaces. For the test variable, we propose local
snapshot spaces. Furthermore, we propose a local spectral
decomposition following our recent work~\cite{ chan2015adaptive},
where we consider minimum energy snapshot functions. We discuss the
discrete stability of the system and its relation to the approximation
properties of the velocity field. The resulting approximation error is
minimized within our multiscale framework by selecting a few
multiscale basis functions. Our numerical results show that we can
stabilize the system using a few test functions for a given trial
space. We describe and analyze several relevant numerical examples
that validate our theoretical results.

\section*{Acknowledgements}
This work is part of the European Union's Horizon 2020 research and
innovation programme of the Marie Sklodowska-Curie grant agreement No
644602. 2014-0191. This publication also was made possible by a
National Priorities Research Program grant  NPRP grant 7-1482-1-278
 from the Qatar National
Research Fund (a member of The Qatar Foundation). The statements made
herein are solely the responsibility of the authors. 
The work described in this paper was partially supported by a grant 
from the Research Grant Council of the Hong Kong Special Administrative Region, China
(Project No. CUHK 14301314).

\bibliographystyle{plain}
\bibliography{references}
\end{document}